\documentclass{article}

\usepackage{graphicx,xcolor}
\usepackage{comment}
\usepackage{multicol}
\usepackage{amsmath,amsfonts}
\usepackage{graphicx}
\usepackage{epstopdf}
\usepackage{algorithmic}
\usepackage{hyperref}
\usepackage[letterpaper,margin=3.00cm]{geometry}                		
\usepackage{authblk}
\usepackage{amsopn}

\newcommand{\R}{\mathbb{R}}

\newcommand{\p}{\partial}

\newcommand{\e}{\mathrm{e}}
\renewcommand{\d}{\mathrm{d}}
\newcommand{\g}{\mathbf{g}}

\newcommand{\I}{\mathbf{I}}

\title{A boundary integral equation method for the complete electrode model in electrical impedance tomography with tests on experimental data}

\author[1]{Teemu Tyni}
\author[2]{Adam R Stinchcombe}
\author[2]{Spyros Alexakis}

\affil[1]{Unit of Applied and Computational Mathematics \newline Faculty of Information Technology and Electrical Engineering, University of Oulu \newline Pentti Kaiteran Katu 1, 90014 Oulu, Finland \newline (\texttt{teemu.tyni@oulu.fi}).}
\affil[2]{Department of Mathematics, University of Toronto \newline 40 St.~George Street, Toronto ON, M5S 2E4, Canada \newline (\texttt{stinch@math.toronto.edu}, \texttt{alexakis@math.toronto.edu}).}

\date{}

\begin{document}

\maketitle

\begin{abstract}
    We develop a boundary integral equation-based  numerical method to solve 
    for the electrostatic potential in two dimensions, 
    inside a medium with piecewise constant conductivity, where
    the boundary condition is given by the complete electrode model (CEM). The CEM is seen as the most accurate model of the physical setting where electrodes are placed on the surface of an electrically conductive body, and currents are injected through the electrodes and the resulting voltages are measured again on these same electrodes. The integral equation formulation is based on expressing the  electrostatic potential as the  solution to a finite number of Laplace equations which are coupled through boundary matching conditions. This allows us to re-express the solution in terms 
    of single layer potentials; the problem is thus re-cast as 
     a system of  integral  equations on  a finite number of smooth curves. 
     We discuss an adaptive method for the solution of the resulting system of mildly singular integral equations. This forward solver is both fast and accurate. We then present a numerical inverse solver for electrical impedance tomography (EIT) which uses our forward solver at its core.  To demonstrate the applicability of our results we test our numerical methods on an open electrical impedance tomography data set provided by the Finnish Inverse Problems Society.
\end{abstract}


\tableofcontents

\section{Introduction}
In this paper we undertake two separate challenges. The first is to produce a fast and accurate 
novel solver for the electrostatic potential in a 2-dimensional conductive body, where the conductivity 
is assumed to be \emph{piecewise constant}, and where the boundary condition is given by current 
injections according to the {complete electrode model}.  This setting is often studied in connection 
with 
the electrical impedance tomography (EIT) inverse problem. The second aim of this paper is to use this 
forward solver to then build a fairly straightforward EIT inverse solver which we test on a set of open 
data for some real-world EIT experiments.

We will start by briefly introducing the setting that we consider for both our forward solver and our inverse solver.

\subsection{The setting for our problem: Body configurations and experiments}

The objects of study in the first part of our paper will be solutions to the electrostatic potential 
equation:
\begin{equation}\label{thePDE}
\nabla\cdot (\sigma(x)\nabla u)=0,
\end{equation}
for $u$ defined over a compact domain $\Omega_0\subset \mathbb{R}^2$ with smooth boundary. The electric conductivity $\sigma(x)$ of the conductive body $\Omega_0$ is assumed to be 
\emph{piecewise constant}.  

In the particular setting we consider here, one can visualize the conductive body $\Omega_0$ as consisting of ``islands'' $\Omega_i$ of different shapes and different conductivities sitting inside a conductive medium. A number of electrodes $L$ are placed at the boundary of the conductive body, and direct or alternating currents $\vec{I}$ are injected precisely at those electrodes. The currents here consist of $L$ numbers $\vec{I}=(I_1,\dots, I_L)$ subject to the restriction that $\sum_{i=1}^LI_i=0$. The voltages (up to a grounded electrode)  are measured again at the $L$ electrodes. So these are again captured by  $L-1$ numbers (since one of the electrodes is grounded) which we list out in  $\vec{V}=(U_1,\dots, U_{L-1})$.

Thus a particular configuration $\vec{B}$, among a set ${\cal B}$ of permitted configurations, of conductive bodies inside a media, with certain fixed electrodes $e_1,\dots e_L$ touching the boundary determines a current-to-voltage map $\operatorname{R}_{\vec{B}}(\cdot)$, which we denote by
\[
\vec{V}=\operatorname{R}_{\vec{B}}(\vec{I}).
\]
The map $\operatorname{R}_{\vec{B}}$ is called the resistance matrix \cite{CEM-exist-unique2}. We use the subscript $\vec{B}$ to stress the dependence of the map on the configuration $\vec{B}\in {\cal B}$. The map is linear in $\vec{I}$, but its dependence on $\vec{B}$ is of course very much non-linear. 

Finding the configuration $\vec{B}$ using the knowledge of $\operatorname{R}_{\vec{B}}$ is the topic of the EIT inverse problem. This question has been studied very extensively both in the pure math literature (see the pioneering works \cite{SylvesterUhlmann,Nachman}, which have led to many generalizations) and also in more real-world settings. We note that the transition from mathematical solutions to real-world reconstructions is often a difficult one, owing to the well-known ill-posedness of the problem. If one allows a very large space of possible conductivities $\sigma(x)$, then it is known~\cite{M-instability, Uhlmann-survey1} that changes of the conductivity can result in comparatively much smaller changes in the measured voltages map $R_{\vec{B}}(\vec{I})$, for a given injection $\vec{I}$. Sensitivity of measured voltages to perturbations of the conductivity $\sigma$ have been considered in \cite{DHMV2021}.

If more restrictions are imposed  on the conductivity inside the body then better stability is possible. In~\cite{AV-lipschitz} it was shown that the recovery of piecewise constant conductivity with finitely many bodies is Lipschitz stable (at least from the full Dirichlet-to-Neumann map). The Lipschitz constant, however, is not estimated. In fact the proper understanding of the stability of reconstruction should also depend on the norms in which (differences in) conductivities are measured. We do not touch on this matter further, but turn to potential uses of EIT.

\subsection{EIT as an imaging tool}

EIT has been proposed and developed as a possible tool in medical imaging over many decades. A resource on possible uses, and the state-of-the art in 2005, is~\cite{Holder}. Let us very briefly describe some areas of more intense inquiry.

As discussed in~\cite{M-instability,Uhlmann-survey1}, the severe ill-posedness of EIT, at least when one seeks to reconstruct conductivities that could a priori lie in a very large dimensional space seems to suggest that one cannot \emph{in general} expect high-resolution reconstructions.  However, its non-invasive and inexpensive nature makes it an attractive candidate for medical imaging. 

Medical imaging applications of EIT include imaging lung function~\cite{Isaacson1, Isaacson2}; studies towards breast cancer detection~\cite{breast-cancer-detection}; and towards brain imaging~\cite{brain-activity-monitor, head-imaging}. A further very interesting use of EIT is when one does not seek full imaging, but rather to detect very coarse features inside a body. One such application is to decide whether a stroke is due to haemorraging or a blood clot--in other words, using EIT towards stroke detection~\cite{stroke-detection}.  

We also mention that the EIT problem can have variants, where one has \emph{less} information than discussed in the previous section. For example, in some real-world instances one can have the shape of the boundary $\partial\Omega_0$ \emph{not known}. (One such instance is when the boundary is changing in time, due to breathing during lung monitoring). Also, the contact impedances of the electrodes at the interfaces with the body can \emph{not} be assumed, and must also be determined. A reference for the reconstruction of these parameters is~\cite{CEM-simulatenous-2021}. 

\subsection{EIT recast as an optimization problem}
We return to the setting where the space of possible conductivities $\cal B$ is assumed to be a finite-dimensional space. The EIT problem in this setting is to determine (to a decent approximation) the configuration of bodies $\vec{B}$ from a sufficient number of samples to the current-to-voltage map $\operatorname{R}_{\vec{B}}$.  In more practical terms, we consider the setting where the space of bodies $\vec{B}$ is of finite dimension (say $D$), and we consider a finite set $\cal I$ of injections $\vec{I}$. For each injection, we sample the voltages at all electrodes, up to grounding.  Thus we have $|{\cal I}|\cdot (L-1)$ pieces of data. Assuming  $|{\cal I}|\cdot (L-1)\gg D$, we can hope to reconstruct the configuration $\vec{B}\in {\cal B}$. 

One way to achieve this is by a sufficiently accurate solver for the forward problem. Indeed, the most ``direct'' method in the inverse problem is to consider the same experiments 
\[
{\cal I}\ni\vec{I}\mapsto \vec{V}=\operatorname{R}_{\vec{B}}(\vec{I}), 
\] 
for $\vec{B}$ in the configuration space ${\cal B}$ and to minimize the error between the ``measured'' current-to-voltage map samples and those that \emph{would} be the samples, among all candidate configurations $\vec{B}\in{\cal B}$. Perhaps the simplest measure of difference between to current-to-voltage maps is a sum of squares between the differences of the voltages measured for the same injection. 

Assuming we had an accurate solver for the problem for all configurations $\vec{B}\in {\cal B}$, and given a \emph{true} configuration $\vec{B}_{\rm true}$, and a fixed set ${\cal I}$ of current injections, we can consider the $L^2$-cost function:
 \[
{\rm Cost}_{\cal I}(\vec{B}, \vec{B}_{\rm true})= \sum_{\vec{I}\in {\cal I}} \|\operatorname{R}_{\vec{B}}(\vec{I})- \operatorname{R}_{\vec{B}_{\rm true}}(\vec{I})\|^2.
 \]
Thus, any solver for the forward problem can yield an  EIT inverse solver: One would seek to identify the configuration $\vec{B}\in{\cal B}$ which minimizes the above cost function, among all $\vec{B}\in {\cal B}$. (Or more generally, some suitable finite-dimensional subspace of $\vec{B}$, if one had a priori information on $\vec{B}_{\rm true}$). To make this method practically applicable, one needs forward solvers that are fast and accurate. Speed is needed since even a modest \emph{iterative} minimization algorithm should be expected to take about 100 steps. Each step might involve finding the gradient of a function of many variables (at least 10, perhaps);  and each evaluation of the function involves solving at least 15 forward problems, even for very simple models with relatively few (i.e.~16) electrodes. High accuracy is also essential, due to the severe ill-posedness of the EIT inverse problem. As discussed, this roughly means that small errors in the measurement lead to large errors in the reconstruction. 

Our second result here is precisely to construct such an inverse solver out of our new forward solver, and to test its performance on experimental data. The  data we will consider is that produced by the Finnish Inverse Problems Society (FIPS) in  \cite{OpenData}.  Of course, measurement error in the data can affect the quality of the found minimum, in view of the ill-posedness of the inverse problem. There seems to be some amount of measurement error in the FIPS data, which we will discuss later. 

\subsection{Alternative approaches to the EIT inverse problem}
As an alternative to our inverse solver (which relies on minimizing a simple sum-of-squares function), one could consider more cleverly chosen cost functions, or seek to incorporate a priori knowledge on the possible solutions by adding a regularizer to the cost function, or by using Bayesian inference. Some references on such approaches are \cite{KaipioSomersalo,Stochastic-PDE}. Further approaches to the inverse problem have been proposed that do not involve minimization of a cost function at all. Indeed, the factorization method~\cite{BH-factorization-eit1,BH-factorization-eit2,BH-factorization-eit3,Kirsch-factorization1,Kirsch-factorization2} (see also factorization method for the CEM~\cite{HHS-factorization-CEM,LHH-factorization-CEM}), monotonicity based methods~\cite{Harrach-monotonicity}, as well as the implementation of the d-bar methods (specific to two dimensions)
take samples of the current-to-voltage map to \emph{directly} reconstruct $\vec{B}$. The factorization method 
is also taylored to locating sharp jumps in conductivity; the d-bar method has also been tried in such experimental settings;  results for those methods are available in~\cite{Isaacson1,Isaacson2,d-bar1,d-bar2}. 

We note that  one could imagine combining such ``direct'' methods with an inverse solver based on minimizing a cost function, by using a direct method \emph{first}, to find an suitable ``first guess'' for the desired $\vec{B}_{\rm true}$, say $\vec{B}_{\rm first-guess}$. This $\vec{B}_{\rm first-guess}$  can then be used as the initial point for a suitable optimizer. 

We also mention recent progress on applying learning techniques to identify better cost functions, or better regularizers \cite{Carola}. It is possible that one can enhance the rudimentary inverse solver we present here, by using ``learned'' data from a forward solver.  

However, our aim in this paper is not to design a very good inverse solver based on, for example, a smarter choice of cost function, or using a more tailored optimization software. Instead, we wish to complement the new (fast and accurate) forward solver we obtain with a rather basic inverse solver, which we test on experimental data.

In \S \ref{sec:CEM}  we introduce the PDE for the electrostatic potential, and 
recast it as a coupled system of integral equations. In \S \ref{sec:BIEM} we present the 
numerical solver that implements the system derived in \S\ref{sec:CEM}, and discuss the adaptive refinements, as well as its convergence and accuracy. In \S\ref{sec:experimental}  we study the 
FIPS open data set, and define the cost function that we will be minimizing for our inverse solver. We then test our inverse solver on most of the FIPS data in 
\S\ref{sec:real-data-reconstructions}.

\section{The CEM PDE for the electrostatic potential recast as a system of integral equations}\label{sec:CEM}

The complete electrode model is widely considered to be the most accurate 
model for electrical impedance tomography, where injections of current are performed through electrodes that touch the boundary of the body. 
The CEM setting has been studied by many authors, both theoretically (e.g.~\cite{CEM-exist-unique1,CEM-exist-unique2}) and also numerically (see for instance~\cite{AC2007,DLA2015,hp-adaptive-cem,FEM-convergence}).

We present only the very basics here in order to then derive the integral equation formulation of this problem.

\subsection{The complete electrode model as a system of coupled Laplace equations}

Let $\Omega_0\subset\mathbb{R}^2$ be a compact domain with smooth boundary. 
We let  $\Omega_i\subset \Omega_0$, $i=1,2,\dots, N$ be domains 
inside $\Omega_0$, which represent ``bodies'' with smooth boundaries.

In the interior of the domain $\Omega_0$, the electrostatic potential $u$ satisfies the elliptic PDE
\begin{equation}\label{the PDE}
\nabla\cdot (\sigma(x)\nabla u)=0,
\end{equation}
in which $\sigma(x)$ is assumed to be a constant $\sigma_i$ inside each of the bodies $\Omega_i$ when they do not overlap, and equal to $\sigma_0$ in $\Omega_0\setminus (\bigcup_{i=1}^N\Omega_i)$. 

The boundary condition on $u$ on $\partial\Omega_0$ is as follows. Away from all electrodes $e_k$, one has zero current crossing into $\Omega_0$, so the condition is
\[
\sigma\partial_n u=0,\quad x\in\p\Omega_0\setminus \cup_{k=1}^L e_k,
\]
where $n$ denotes the outward unit normal vector of $\Omega_0$. On the electrodes, the condition is two-fold. The function should have constant Dirichlet data if the contact impedance of an electrode is zero. In general one expects a (small) contact 
impedance $z_k>0$ on each electrode $e_k$. So the condition of constant voltage on the electrode $e_k$ reads
\[
u(x)+\sigma z_k\partial_n u(x)=U_k,\quad x\in e_k,
\]
in which $U_k$ is an unknown constant (voltage) on the perfectly conducting electrode.
We then prescribe the total current through the electrodes
\[\int_{e_k} \sigma \p_n u~\d l_x = I_k,
\]
with the only constraint that the total current through $\partial \Omega_0$ should vanish, $\sum_{i=1}^L I_k=0$. As discussed, the current-to-voltage map\footnote{We include the zero voltage on the one grounded electrode.} is then the map 
\[
R:(I_1,\dots I_L)\mapsto (U_1,\dots U_L). 
\]

\begin{figure}
    \centering
    \includegraphics[width=0.5\textwidth]{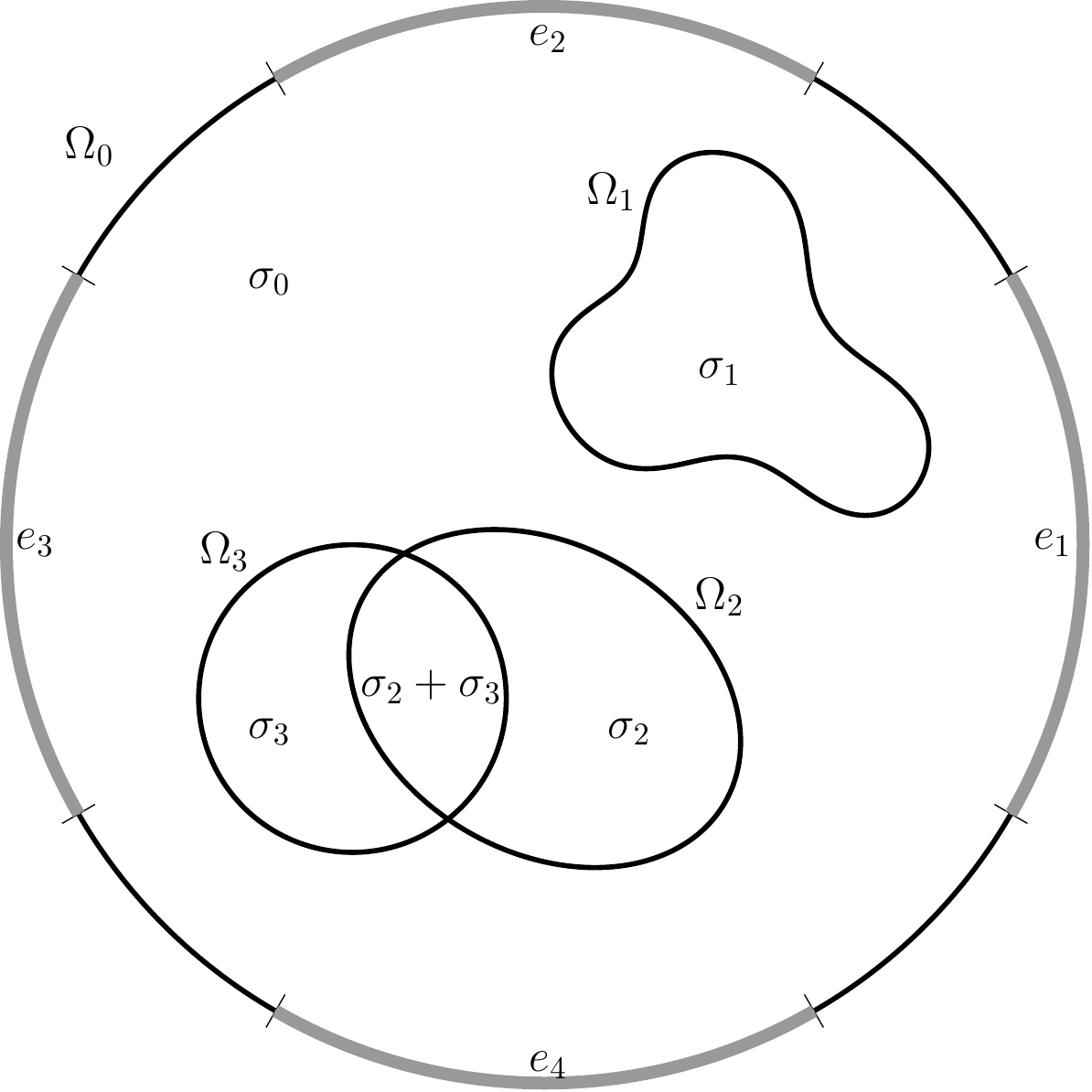}
    \caption{An example configuration of a domain $\Omega_0$ with three bodies $\Omega_1,\Omega_2$ and $\Omega_3$ with their associated conductivities $\sigma_0,\ldots,\sigma_3$, respectively. On the intersection $\Omega_2\cap\Omega_3$ we define the conductivity to be $\sigma_2+\sigma_3$. In this example there are four electrodes $e_i$, $i=1,2,3,4$, whose locations on $\p\Omega_0$ are marked with gray color.\label{fig:domainplot}}
\end{figure}

It is useful to write an equivalent formulation of the above problem, recasting the equation as a set of harmonic functions away from all interfaces $\partial\Omega_i$, with suitable matching condition at those interfaces. The complete electrode model (CEM) problem  for \eqref{the PDE} with $N$ domains $\Omega_i$ inside a domain $\Omega_0$ (see Figure~\ref{fig:domainplot}) is given by the following system:
\begin{equation}\label{eq:complete-electrode-model}
\begin{cases}
    \Delta u = 0, & x\in\Omega_0\setminus \cup_{i=1}^{N} \p\Omega_i,\\
    \sigma \p_n u = 0, & x\in\p\Omega_0\setminus \cup_{k=1}^L e_k,\\
    u+\sigma z_k \p_n u = U_k, & x\in e_k,\, k=1,\ldots,L,\\
    \int_{e_k} \sigma \p_n u~\d l_x = I_k,& k=1,\ldots,L,\\
    [u] = 0,&\text{on }\p\Omega_i,\, i=1,\ldots,N,\\
    [\sigma\p_n u] = 0,& x\in\p\Omega_i,\, i=1,\ldots,N.
    \end{cases}
\end{equation}
We explain the notation used above: 
To each body $\Omega_i$, $i=0,1,\ldots,N$, we assign a conductivity $\sigma_i> 0$. To allow the bodies to overlap, we use the convention 
$$
\sigma(x) := 
\begin{cases}
\sum_{i=1}^N \sigma_i\chi_i(x),& x\in \cup_{i=1}^N\Omega_i\\
\sigma_0, & x\in \Omega_0\setminus \cup_{i=1}^N\Omega_i,
\end{cases}
$$
in which $\chi_i(x)=1$ if $x\in \Omega_i$ and $\chi_i(x)=0$ otherwise is an indicator function for the body $i$.
We use the notation $[\xi](x) = \xi^+(x)-\xi(x)^-$, where $\xi(x_0)^+$ (resp. $\xi(x_0)^-$) is the limit of $\xi(x)$ as $x\to x_0$ from $\Omega_0\setminus \Omega_i$ (resp. from $\Omega_i$). We assume that the interfaces $\p\Omega_i$, $i=1,\ldots,N$, intersect transversely and at most at finitely many points; the interface condition --- the last equation of \eqref{eq:complete-electrode-model} --- is well-defined everywhere apart from these possible intersection points.

The electrode potentials $U_k$ in \eqref{eq:complete-electrode-model} are unknown constants that are implicitly determined by the known electrode currents $I_k$. Each electrode has a known, constant, positive contact impedance $z_k$. Due to charge conservation and so that the problem will have a solution,  we require that there is no net current into the domain, i.e. 
\[
\sum_{k=1}^L I_k = 0.
\]
Finally, for uniqueness of the solution we impose $\sum_{k=1}^L U_k = 0$.
Under these conditions, well-posedness of \eqref{eq:complete-electrode-model} in $H^1(\Omega_0)$ is well known, see for example \cite{CEM-exist-unique1,CEM-exist-unique2}.

\subsection{Derivation of the boundary integral equations}

As discussed above, the ``real world'' EIT problem occurs in spatial dimension 3. It is possible however to produce experimental data that corresponds to spatial dimension two when all objects involved (the bodies, the medium in which they are placed, as well as the electrodes) are defined over a cylindrical region $\Omega\subset\mathbb{R}^3$, with $\Omega= D\times [0,h]$, $D\subset \mathbb{R}^2$,  and all objects enjoy cylindrical symmetry: i.e.~they are invariant under changes in the coordinate $x_3\in [0,h]$, and moreover the boundary conditions on the currents $\vec{J}$ at the two boundaries $\{x_3=0\},\{x_3= d\}$ are $\vec{J}=\vec{0}$. This setting can be achieved experimentally (see the data set below), and is the 3D ``lift'' of the 2D setting that we consider. In particular our forward solver applies only to 2D. Our approach to this problem is via single-layer potentials, which we present next. 

We note that an earlier paper by two of the authors, joint with K.~Bower and K.~Serkh presented a boundary integral equation method for piecewise constant conductivities in 2D, where a more ``classical'' Neumann boundary condition is imposed~\cite{Bower}. A similar approach of using single layer potentials to solve a Robin boundary value problem for the Laplace equation was taken in~\cite{BIEM-robin}. Beyond the fact that the CEM is a realistic model in the presence of electrodes, we add that the numerics need to be more attuned to the presence of singularities of the solution where the electrodes terminate. Furthermore, in this paper we allow for the bodies $\Omega_i$ in the interior to overlap, which requires novelty in the numerics, to obtain precision near the intersection points of the boundaries $\partial\Omega_i, \partial\Omega_j$. This particular idea was first developed in~\cite{Bower-thesis}. 

We represent the solution $u$ of \eqref{eq:complete-electrode-model} as a sum of single layer potentials
\begin{equation}\label{eq:solution u sum}
 u(x) = \sum_{j=0}^N S_{\p\Omega_j}[\gamma_j](x),
\end{equation}
with the single layer potential $S_{\p\Omega_j}[\gamma_j]$ defined via
\[
S_{\p\Omega_j}[\gamma_j](x):= \int_{\p\Omega_j} G(x,y)\gamma_j(y)~\d l_y,
\]
with the weakly singular kernel
\[
G(x,y) := \frac{1}{2\pi} \log(|x-y|).
\]
We note that the single layer potential is a harmonic function, $\Delta S_{\p\Omega_j}[\gamma_j]=0$. Using the classical jump relations for single layer potentials, we see that
\[
\gamma_i = [\p_n u] = \p_{n^+}u - \p_{n^-}u,\quad \text{on } \p \Omega_i,\, i=1,\ldots,N.
\]
The normal derivative of the single layer potential is given 
 for $x\in\p\Omega_i$ by
\[
\p_{n^\pm} S_{\p\Omega_i}[\gamma_i](x) = K^*_{\p\Omega_i}[\gamma_i](x) \pm \frac{1}{2}\gamma_i(x)
\]
(see e.g., \cite{BIEM-book1,Kress-book1}), which we may write for $x\in\p\Omega_i$, $i=0,\ldots,N$ as
\[
\p_{n^\pm} S_{\p\Omega_j}[\gamma_j](x) = \int_{\p\Omega_j} K(x,y) \gamma_j(y)~\d l_y 
\pm \gamma_j(x)\delta_{ij},
\]
in which $\delta_{ij}$ is the Kronecker delta symbol.
Above $K^*$ is the adjoint of the Neumann-Poincar\'e operator given by
\[
K^*_{\p\Omega_i} := \frac{1}{2\pi} \int_{\p\Omega_i}
\frac{(x-y)\cdot n(x)}{|x-y|^2} \gamma_i(y)~\d l_y,\quad x\in \p\Omega_i.
\]
Let
% The factor of 1/2\pi was missing. It has been added on 22 Aug 2023
\[
K(x,y) := \frac{(x-y)\cdot n(x)}{2\pi|x-y|^2}.
\]
The kernel $K(x,y)$ has a removable singularity on the diagonal, provided $\p\Omega_i$ is $C^2$-smooth. By defining $K(x,x) = \frac{\kappa(x)}{4\pi}$, where $\kappa(x)$ is the curvature at point $x\in\p\Omega_i$, we find that $K$ is a continuous function.

Substituting~\eqref{eq:solution u sum} into~\eqref{eq:complete-electrode-model}, we obtain the integral equation formulation of the CEM.\\
\textbf{Off electrodes}:
\begin{equation}\label{eq:CEM off electrodes}
    -\frac{1}{2}\sigma(x) \gamma_0(x)
    +
    \sigma(x) \sum_{j=0}^N \int_{\p\Omega_j} \big(K(x,y)-\delta_{1j}\big) \gamma_j(y)~\d l_y = 0,\quad x\in \Omega_0\setminus\bigcup_{k=1}^L e_k.
\end{equation}
\textbf{On electrodes}:
\begin{equation}\label{eq:CEM on electrodes}
    -\frac{1}{2}z_k\sigma(x) \gamma_0(x)
    + 
     \sum_{j=0}^N \int_{\p\Omega_j}\Big( G(x,y) + z_k\sigma(x) K(x,y)\Big) \gamma_j(y)~\d l_y = U_k,
\end{equation}
for $x\in  e_k,\, k=1,\ldots,L$;
\begin{equation}\label{eq:CEM integrals of neumann}
     \int_{e_k} \Bigg(-\frac{1}{2}\sigma(x) \gamma_0(x)
    +
    \sigma(x) \sum_{j=0}^N \int_{\p\Omega_j} K(x,y) \gamma_j(y)~\d l_y\Bigg) \d l_x%\\
    %&=\int_{e_k} \sigma \p_n u ~\d l_x 
    = I_k,
\end{equation}
for $k=1,\ldots,L$.\\
\noindent
\textbf{On interfaces}:
\begin{equation}\label{eq:CEM inner interfaces}
    \frac{1}{2}(\sigma^+(x)+\sigma^-(x))\gamma_i(x)
    +
    (\sigma^+(x) - \sigma^-(x)) \sum_{j=0}^N \int_{\p\Omega_j} K(x,y) \gamma_j(y)~\d l_y = 0, 
\end{equation}
for $x\in\p\Omega_i,\, i=1,\ldots,N$. In~\eqref{eq:CEM off electrodes}, we replace the kernel $K(x,y)$ by $K(x,y)-\delta_{1j}$. This enforces~\cite{Bower} that the total charge on the outer interface satisfies
\[
\int_{\p\Omega_0\setminus \cup_{i=1}^Le_i} \gamma(x)~\d l_x = 0.
\]
\begin{comment}
\textbf{On inner interfaces}:
\begin{equation}%\label{eq:CEM inner interfaces}
    \frac{1}{2}(\sigma_{p_i}+\sigma_i)\gamma_i(x)
    +
    (\sigma_{p_i}-\sigma_i) \sum_{j=0}^N \int_{\p\Omega_j} K(x,y) \gamma_j(y)~\d l_y = 0, 
\end{equation}
where $x\in\p\Omega_i,\, i=1,\ldots,N$.
\end{comment}

\section{Numerical implementation of the BIEM}\label{sec:BIEM}
\subsection{Discretization of the integral equations}

The integral equations \eqref{eq:CEM off electrodes} -- \eqref{eq:CEM inner interfaces} are discretized using the Nystr\"om method. We first parametrize the interfaces $\p\Omega_i$ by $q\in[0,1)$ and then discretize these curves by taking a grid of $M_i>0$ points $q^i_j$, $j=0,\ldots,M_i-1$. For accuracy, we use two different discretization schemes on the interfaces, panelled and uniform. Specifically, we split the outer interface into panels ensuring that there are panel breaks at the electrode boundaries. We also split interfaces that are within a specified distance to another interface into panels. Moreover, in the case that two or more interfaces overlap, we break panels at the intersection points that are found with a quasi-Newton method. We subdivide panels neighboring intersections points and electrode boundaries dyadically following~\cite{HELSING20088820}. On panelled interfaces, we use Gauss-Legendre quadrature on each panel. On uniform interfaces, we use trapezoidal rule quadrature for the periodic $q\in[0,1)$. This is exponentially convergent provided the interface is distant from other interfaces~\cite{TrefethenWeideman2014}.

The kernel $G(x,y) + z_k\sigma(x)K(x,y)$ in~\eqref{eq:CEM on electrodes} has a logarithmic singularity and so we use auxiliary nodes~\cite{Hao2014} to accurately evaluate these integrals. Specifically, for each node $q^i_j\in e_k$ on an electrode, we introduce an auxiliary set of nodes $\{x_i\} \in e_k$ on each side of $q^i_j$. These nodes and their associated weights $\{w_i\}$ are a generalized Gaussian quadrature~\cite{bremer2010GGQ} designed to integrate 
$$ \int_0^{\Delta q} f_1(q) \log q + f_2(q)~\d q$$ 
for smooth functions $f_1$ and $f_2$\footnote{This quadrature was computed with \texttt{GGQ} available at \url{https://github.com/JamesCBremerJr/GGQ}.}. This quadrature requires samples of the integrand at the nodes $\{x_i\}$ rather than samples of $f_1$ and $f_2$.

The use of this generalized Gaussian quadrature is enabled by properties of the logarithm. The quadrature does not depend on ${\Delta q}$ since $\log c q = \log c + \log q$ so that ${\Delta q}$ can be scaled away without changing the form of the integrand. Additionally, for $F(q) = | x(q) - x(q^i_j) |$, the singular part of the kernel can be written as 
\begin{equation*}
\log F(q) = \log \left( F(q^i_j) + (q-q^i_j) F'(q_\ast) \right) = \log(q-q^i_j) + \log( F'(q_\ast) )
\end{equation*}
using the mean value theorem and $F(q^i_j)=0$. Thus, our integrands can be written as $f_1(q) \log q + f_2(q)$ for smooth $f_1$ and $f_2$. We assume that sufficient dyadic subdivisions have been made so  that the charge density $\gamma_i$ is well-approximated by a polynomial near the node $q^i_j$. To obtain the values of the charge densities at the nodes $\{x_i\}$ we use Lagrange interpolation, which can be augmented into the system of equations for CEM since this interpolation is a linear operation. This interpolation is well-conditioned since the interpolation nodes are either the Gauss-Legendre nodes on panels or uniform on a periodic $q\in[0,1)$~\cite{Trefethen2019Approximation}.

After discretization, the equations \eqref{eq:CEM off electrodes} -- \eqref{eq:CEM inner interfaces} become a linear system,
\begin{equation}\label{eq:CEM discretized}
A\g = \I,\quad \g\in \R^{M+L},
\end{equation}
in which $M=\sum_{i=0}^N M_i$ is the total number of grid nodes. The equations \eqref{eq:CEM off electrodes}, \eqref{eq:CEM on electrodes}, and \eqref{eq:CEM inner interfaces} yield one block of the matrix $A$. The equations \eqref{eq:CEM integrals of neumann} are then appended to the system as a wide block matrix. Similarly, since the electrode potentials $U_k$ are unknown, they are appended to the system as a tall block matrix. Roughly speaking, the matrix $A$ has the form
$$
A=
\begin{pmatrix}
\mathbf{K} & \mathbf{U}\\
\mathbf{N} & \overline{0}
\end{pmatrix},
$$
in which $\mathbf{K}\in\R^{M\times M}$ contains the (discretized) equations \eqref{eq:CEM off electrodes}, \eqref{eq:CEM on electrodes} and \eqref{eq:CEM inner interfaces}; $\mathbf{U}\in\R^{M\times L}$ contains zeros on rows corresponding to nodes off electrodes and $-1$'s on rows corresponding to electrodes; and $\mathbf{N} \in \R^{L\times M}$ contains the (discretized) equations \eqref{eq:CEM integrals of neumann}. Here $\overline{0}$ is an $L\times L$ zero matrix and the right-hand side is given by the injected currents by $\I=(0,\ldots,0,I_1,\ldots,I_L)^T\in \R^{M+L}$. The matrix $\mathbf{K}$ is dense and by far the largest block of $A$ making $A$ a dense matrix.

Solving for $\g$ in \eqref{eq:CEM discretized} yields directly the values of the voltages $U_k$, $k=1,\ldots,L$, on the electrodes $e_k$. This is useful for our electrical impedance tomography inverse problem since those voltages are all that is needed to evaluate the objective function. It is not necessary to evaluate the solution $u$ of \eqref{eq:complete-electrode-model} at all, in stark contrast to the case of  Dirichlet or Neumann boundary conditions. Additionally, in practice one performs many injections, that is, solves the system $A\g=\I$ for several different right-hand sides $\I$. We solve $A\g=\overline{\I}$ for a \emph{matrix} $\overline{\I} :=(\I_1,\ldots,\I_m)$ with Gaussian elimination, solving the multiple injections simultaneously.

\subsection{Adaptive grid refinement}

To accurately resolve the charge densities and have a robust forward solver, we use adaptive grid refinement. We use the spectral resolution of the charge density on an interface to determine whether or not to increase the number of nodes on that interface.

For a uniform interface with $M_i$ nodes, we expand the charge density $\gamma_i$ in a Fourier sum
$$
\gamma_i(q) = \sum_{k=-K}^{K} c_k e^{2\pi \mathrm{i} k q},\quad c_k = \int_0^1 \gamma_i(q) \e^{-2\pi \mathrm{i} kq} ~\d q
$$
in which $M_i = 2K+1$. If the two highest Fourier coefficients $|c_{K-1}|$ and $|c_K|$ are below a specified threshold, then we say that the charge density has been resolved accurately and $M_i$ nodes are sufficient to resolve the charge density. If the coefficients are not below the specified threshold, we increase the number of nodes and re-solve the linear system until the coefficients $|c_{K-1}|$ and $|c_K|$ are small enough.

On interfaces with panels and composite Gauss-Legendre quadrature, we use the Legendre series on each panel separately. If the last two Legendre coefficients of the charge density are above a specified threshold, we refine the panel. If the number of nodes on a panel exceeds a specified maximum number of nodes we instead split the panel into four smaller panels. This procedure is repeated and the linear system is re-solved until the highest order coefficients of the Legendre polynomials are below a specified threshold on each panel.

\begin{figure}
    \centering
    \includegraphics[width=\textwidth]{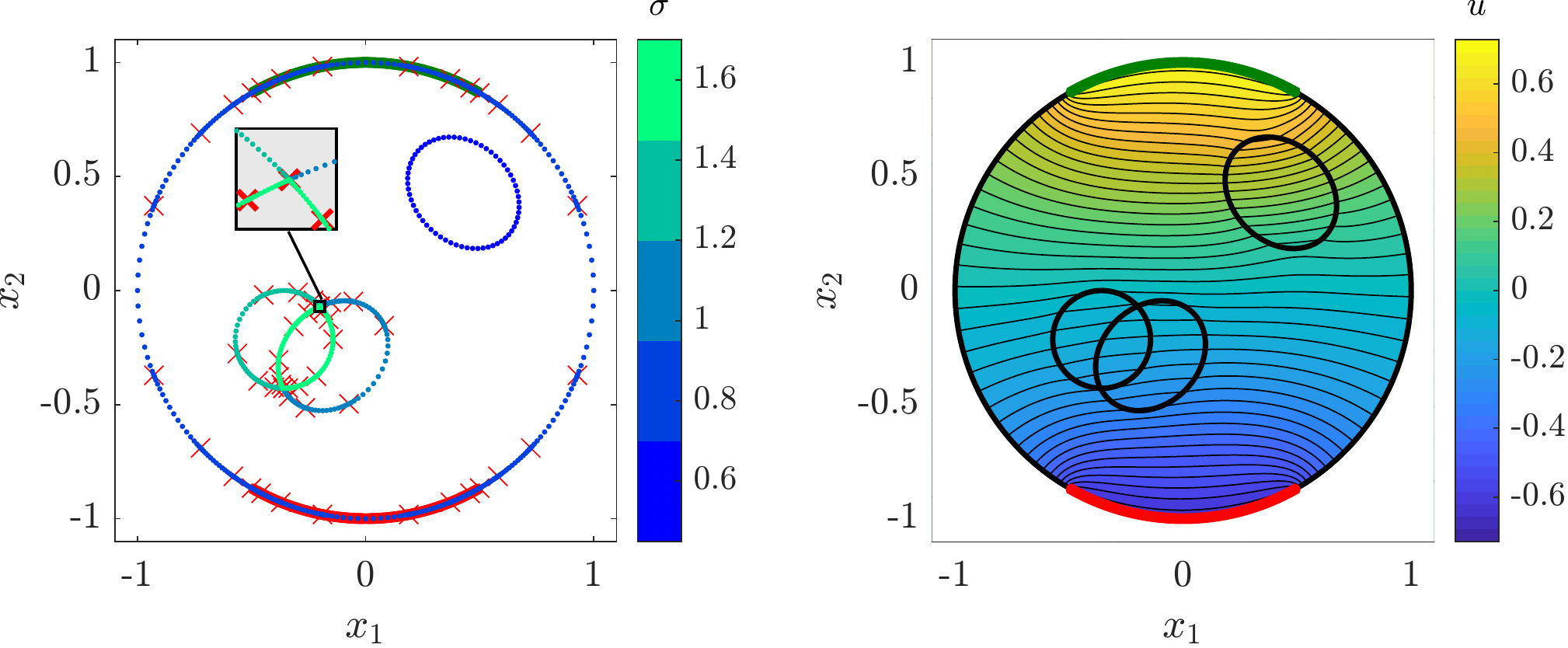}
    \caption{Left: The adaptively selected quadrature grid. The quadrature nodes on interfaces are denoted by dots colored by the conductivity inside the interface on which they are located. The electrode locations are marked with solid red and green curves at the top and bottom of the outer interface. The red crosses mark the panel breaks. Note how our adaptive method refined locations where the interfaces are close to each other.  Right: The computed solution $u$ of \eqref{eq:complete-electrode-model}. We observe that the potential $u$ is approximately constant near the electrodes. Recall that on the overlap region of two bodies the conductivity is defined to be the sum of the conductivities of the bodies, which can be seen in the left figure on the overlap of the two bottom-left bodies. The top-right body has a lower conductivity than the two bodies on the bottom-left.\label{fig:overlappingCurves}}
\end{figure}

Figure~\ref{fig:overlappingCurves} demonstrates adaptive refinement of the 
grid. 

\subsection{Numerical validation and experiments}

\begin{figure}
    \centering
    \includegraphics[width=\textwidth]{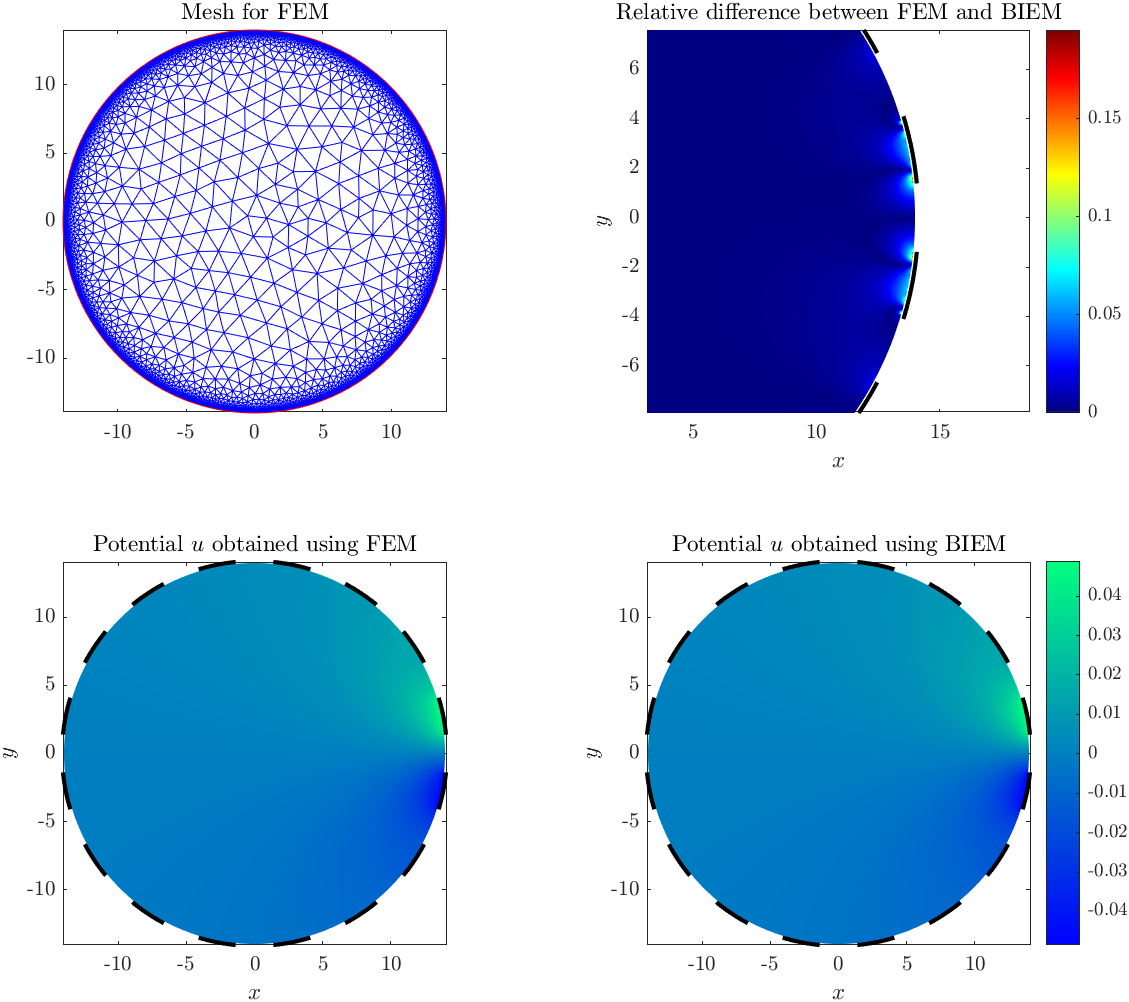}
\caption{Numerical validation of the BIEM by comparison of the potentials $u$ obtained using FEM and BIEM in a simple configuration. The $16$ electrode locations are marked with black lines: active electrodes are located at the right side of the domain. The largest relative difference $\vert u_\mathrm{BIEM}-u_\mathrm{FEM}\vert/\Vert u_\mathrm{FEM}\Vert_{L^\infty(\Omega_0)}$ between the two solutions is attained at the edges of the active electrodes.
\label{fig:femvsbiem}}
\end{figure}

To validate our BIEM solver for CEM, we compared with a simple finite element method (FEM) solver for the CEM. Since we did not attempt to optimize the FEM solver, we do not dwell on a speed comparison of the two methods. However, the BIEM solves a simple problem in approximately 1 to 10 milliseconds while our FEM implementation takes several seconds or minutes. We compare the two approximate solutions $u$ of~\eqref{eq:complete-electrode-model} obtained by FEM and BIEM in Figure~\ref{fig:femvsbiem}, 
we expect to see largest difference in the evaluated solution near the boundary. Despite the close evaluation problem, voltages are accurately computed by the BIEM.

\begin{figure}
    \centering
    \includegraphics[width=\textwidth]{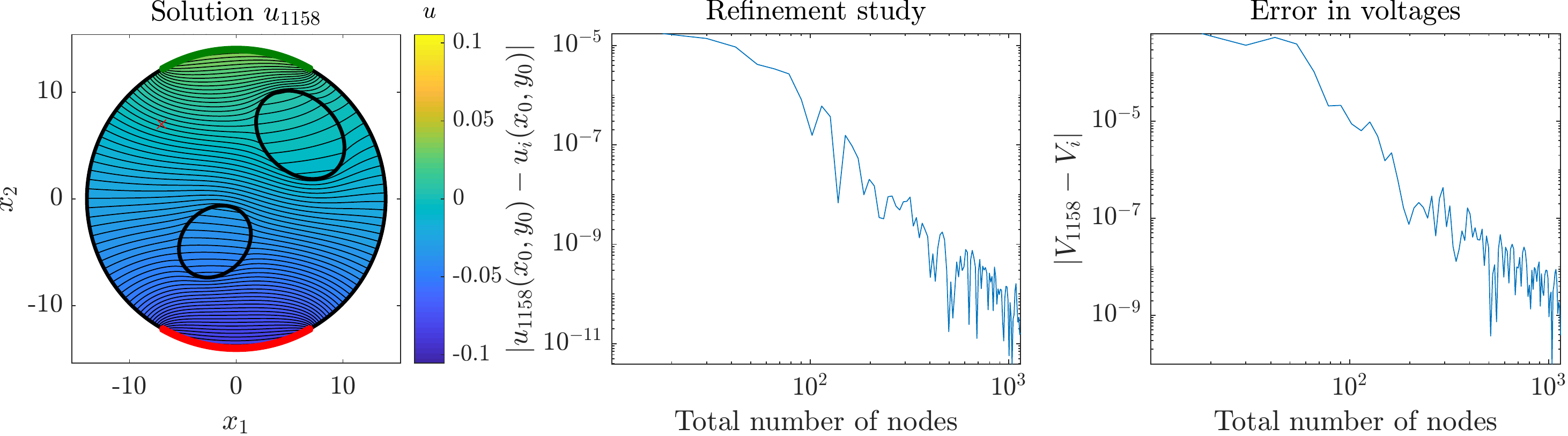}
\caption{Left: the numerical solution $u$ to \eqref{eq:complete-electrode-model} with $1158$ nodes in total. The point $(x_0,y_0)$ is marked  by a red cross. Middle: the refinement study of the absolute differences $\vert u_{1158}(x_0,y_0) - u_{i}(x_0,y_0)\vert$ as a function of the total number of nodes at a point $(x_0,y_0)\in\Omega_0$. Right: the difference $|{\rm V}_{1158}-{\rm V}_{i}|$ between voltages as a function of the total number of nodes. 
\label{fig:errorhistory}}
\end{figure}

To test the convergence of the BIEM we toggled adaptivity off and instead used a fixed number of nodes per panel and per uniform interface. Let $u_i$ denote the numerically obtained solution with $i$ nodes per panel. Figure~\ref{fig:errorhistory} displays the solution to \eqref{eq:complete-electrode-model} in a configuration with two conductive bodies and two electrodes. To study convergence, we calculated the errors between the solutions $u_i$ with successively increasing number of nodes per panel,
and the arising voltages ${\rm V}_i = U_1^{(i)}-U_2^{(i)}$ between the two electrodes. We note that with our adaptive method toggled on we can reach desired accuracy in voltages automatically without the need for an excessively large number of nodes.  

\newpage
\section{Experimental setup and numerical results for the inverse solver}\label{sec:experimental}

In this second part of the paper, we wish to discuss how our new forward solver can be used for solving the EIT inverse problem and to present the results of a rather rudimentary inverse solver on a set of open data~\cite{OpenData} produced by the Finnish Inverse Problems Society\footnote{\url{https://www.fips.fi/}} (FIPS). We start by discussing the experimental setup used to produce the FIPS data set.  

\subsection{The open EIT data set~\cite{OpenData}: Experimental setup}

The experiments were conducted using a tank of circular cylinder shape. The diameter of the
tank was 28\,cm. Sixteen rectangular electrodes (height 7\,cm, width 2.5\,cm) made of stainless steel were attached equidistantly on the inner surface of the tank. The electrodes were marked with orange tapes on the tank wall.

The tank was filled with \textbf{saline} up to the height of 7\,cm, i.e., to the top level of the electrodes.
The measured value of the conductivity was 300\,$\mu$S/cm and temperature was $19^\circ$C. The amplitudes of the injected currents were 2\,mA and their frequencies were 1\,kHz. Corresponding to each current injection, voltages were measured between all adjacent electrodes. Various different phantoms (metal cylinders and plastic inclusions) were placed in the tank in different configurations and the corresponding voltages then recorded.

The voltage measurements are always performed at all pairs of adjacent electrodes.
However, there are different collections of injections: Injections are always between two electrodes; these are performed according to five different patters; in the ``adjacent'' subset 
only adjacent electrodes are excited. In the ``skip $k$'', $k=1,2,3$ subsets only two electrodes separated by $k$ non-excited electrodes are excited; in the ``all against 1'' data set, electrode 1 carries negative current and positive current is injected into each of the other 15.
We note that each of these 5 data sub-sets forms a basis for the set of possible EIT data within the CEM model.

The part of the data set that we use here is \emph{only} that of ``adjacent'' injections of current. Let us number the 16 electrodes by $e_1,\dots, e_{16}$. For convenience, we introduce the convention that $e_{17}=e_1$. 
Then currents $+J$ and $-J$ are injected across adjacent electrodes $e_i, e_{i+1}, i=1,\dots, 16$. 
There are thus 16 injections or ``experiments'' $\vec{I}$, yet only 15 of them are independent (the last one can be reproduced by summing the first 15). For each $i\in \{1,\dots, 16\}$, the voltage 
is measured at all adjacent electrodes $(e_j, e_{j+1})$: For each such injection $\vec{I}$, 
we evaluate the voltage differences 
${\rm V}(j)=U_{j+1}-U_j$, $j=1,\dots 16$. Note that these  voltages are constrained by the condition 
$\sum_{j=1}^{16} {\rm V}(j)=0$. Thus  our data set consists of $15\times 15$ independent elements. 
\medskip

\subsection{The configuration space for the bodies and conductivities, and the inverse solver}
   The inverse solver has at its core the forward solver presented in the first part of
   this paper. 
   It then relies on  a simple minimization (optimization) of a sum-of-squares function ${\rm Cost}_{\cal I}(\vec{B}, \vec{\sigma})$
   which we now define in detail.
   
   Firstly, the configuration space ${\cal B}$ of $N$ 
   piecewise constant conductivities $\vec{B}$  
   will be written as the product of $N$ ``bodies'' $\vec{B}_{\rm shapes}$ and ``conductivity inside the bodies'', $\vec{\sigma}$. 
   So in particular:
   \[
\vec{B}=(\vec{B}_{\rm shapes},\vec{\sigma}),
   \]
where $\vec{B}_{\rm shapes}\in{\bf B}$, and ${\bf B}$ stands for the ``configuration space'' of $N$ ``bodies'' inside a disk; $\vec{\sigma}\in \mathbb{R}^{N+1}_+$ is the parameter space of the conductivities of each of the $N$ bodies, and $\sigma_{N+1}$ is the conductivity of the medium (saline in the case of these experiments).
       
The ``configuration space'' of bodies ${\bf B}$ in this particular situation consists of $A$ bodies which are ellipses, and $B$ bodies which are polygons; in particular either triangles or convex quadrilaterals. 
(We take $A,B$ as known below). We note that we allow different bodies to overlap. In the case of overlap, the conductivity over the overlap region is taken to be the sum of the conductivities assigned to the bodies. 

Each configuration $(\vec{B}_{\rm shapes},\vec{\sigma})$ would then yield a set of $15\times 16$ voltages, if we were to perform the same $16$ injections of current on adjacent electrodes and measured the same $16$ voltages differences as just described: Letting ${\cal I}_i, i=1,2,\dots, 16$ be the 16 injections and $j\in \{1,\dots, 15\}$ we  define 
\[
{\rm V}(j, I_i, \vec{B}_{\rm shapes}, \vec{\sigma})
\]
to be the voltage difference we \emph{would} measure between $e_j, e_{j+1}$ if we performed the injection of current $I_i$ on the configuration of bodies given by $(\vec{B}_{\rm shapes},\vec{\sigma})$.
    
We let 
\[
{\rm V}^{\rm experimental}(j, I_i)
\]
be the \emph{experimentally measured} voltage difference between $e_j, e_{j+1}$ for the injection of current $I_i$. 

The very \emph{direct} approach to the inverse problem that we then take is to  minimize the function
     \begin{equation}\label{cost final}
{\rm Cost}[(\vec{B}_{\rm shapes},\vec{\sigma})]:= 10^6\cdot\sum_{i=1}^{16}\sum_{j=1}^{15}[
{\rm V}(j, I_i, \vec{B}_{\rm shapes}, \vec{\sigma})
-{\rm V}^{\rm experimental}(j, I_i)]^2
     \end{equation}
 over the configuration space $(\vec{B}_{\rm shapes},\vec{\sigma)}\in {\bf B}\times\mathbb{R}^{N+1}_+$. 
 (We will also be minimizing some variants of this cost function, obtained by a suitable 
 \emph{pre-conditioning} of the measured voltages $V^{\rm experimental}(j,I_i)$--see \eqref{cost corrected} below). 

The minimization is performed using optimization software; we use \textsc{Matlab}'s \textsc{fmincon}-function. 
This involves certain choices for the \textsc{fmincon} function, most importantly an initial guess and constraints on the bodies; but also running the inverse solver involves the specifying of parameters for the forward solver that will be used within the inverse solver. Scaling the cost by a large number (the $10^6$ in \eqref{cost final}) improves the effectiveness of the optimization software.

We describe our choices for the \textsc{fmincon} function in brief. Firstly, we choose to \emph{freeze} some of the parameters and not solve for them. In the results we present, we are trying to ``find'' $N$ bodies, each of which is  either an almost perfect conductor (the metal rings) or complete insulators (plastic polygons--either triangles or squares); so in most instances we just \emph{set} the conductivities $\sigma_i$ corresponding to the almost-perfect-conductors/almost-insulators to very large/very small \emph{artificially chosen} numbers, and minimize the cost function ${\rm Cost}[(\vec{B}_{\rm shapes},\vec{\sigma})]$, when those parameters are held \emph{fixed}. We found that this minimization yields almost the same results for $\vec{B}_{\rm shapes}$, compared to when $\vec{\sigma}$ is allowed to vary.

The only constraint we impose is that the bodies should stay inside the (fixed) disk $\Omega_0$, and that $\sigma_i\ge 0$. In particular, during steps in the iteration we can have two different bodies overlapping. This is not  a problem, since in our forward solver works seamlessly in this setting. (Recall we make an auxiliary choice that the conductivity over the overlap region should be the sum of the conductivities of the overlapping bodies). This is a powerful feature of our forward solver, since allowing bodies to overlap allows the inverse solver to take steps where different bodies  move through each other and eventually separate again.  

The space of parameters when we have $A$ ellipses and $B$ triangles, as well as $(N+1)$ unknown conductivities is $[5\times A+6\times B+(N+1)]$-dimensional. If we have $B$ quadrilaterals instead of $B$ triangles it is $[5\times A+8\times B+(N+1)]$-dimensional. The last $N+1$ dimensions correspond to the elements of $\vec{\sigma}\in \mathbb{R}^{N+1}_+$, when all of them are variable.

Here the benefits of our forward solver are evident: One solution of the forward solver takes  0.01 seconds at most, which allows the possibility to converge to an approximate (local) minimum in a very short time. Since we just wish to give a report on the performance of the inverse solver within short time constraints, we only allow the inverse solver 
to run for 128 steps.  We also constrain the initial condition $\vec{B}_{\rm inital}$ to be ``somewhat'' close to the correct position, at least visually. We note that the optimizer at times takes large steps in the parameter space $\cal B$, so its intermediate steps can move much farther from the ``visually correct'' values compared to the initial guess. 

We present the results of the inverse solver for a large part for the experiments published in the FIPS open data set in the next section. Before presenting these results, we make a few remarks.

\subsection{The issue of inaccurate measurements}  As in all real world experiments, the measurements should be expected to involve some error. In the FIPS data set this is apparent in two ways described below. Recall that each subset of measurements (``adjacent'', ``skip $k$'', $k=1,2,3$ and ``all against 1'') form a basis for the space of all possible measurements. 
In particular, had the measurements had perfect precision, \emph{any one} of these data subsets would suffice to reconstruct the others, by simple linear algebra. This fact enables us to detect some errors in the measurements.

Firstly, 
consider separately each of the five 
lists of experiments in the FIPS; we can call each of these a data ``subset''; they are encoded in a $16\times 16$ matrix. Each of these data subsets \emph{viewed on its own} is subject to the constraints that each of the  corresponding rows and columns should sum to zero. This reflects the two properties that given one fixed injection between any pair $e_a, e_{a+1}$, then the corresponding voltages should satisfy $\sum_{i=1}^{16}[U_{i+1}-U_i]=0$ (taking $U_{17}=U_1$). And also that summing over any row calculates the voltage difference when the  injection is zero. Nonetheless, the latter property of elements of rows summing to zero does \emph{not} hold in any of the data subsets. So this is a first evidence of errors in measurement within each of the five data subsets. The second evidence comes by comparing \emph{different} data subsets. 

Since each of the five data subsets is a basis for all possible experiments, one should in  particular 
obtain from it all the experiments in all the other  four data subsets. However, the data obtained, say, from the ``adjacent'' data subset fails to reproduce the experiments in the remaining data sets. We note that the mismatch between the measurements that should have agreed can at times be significant: seeking to reproduce the skip-1 data set using the adjacent data set as a basis can lead to approximately $0.02\,V$ ($1\,\%)$ voltage difference between the different data sets, when the maximum sought voltages are approx. $1.8\,\mathrm{V}$. % For example inject at electrode 5 in phantom1.2
\medskip

\emph{Data subset used:} In the results we present below, we rely \emph{exclusively} of the 
``adjacent'' data subset. In particular the cost function \eqref{cost final} is constructed using the 
experimental Voltages $V^{\rm experimental}(j,I_i)$ corresponding to the adjacent data set. 
However, to illustrate the effect of measurement errors on our reconstruction we will construct 
\emph{two different} cost functions and seek to minimize these two functions \emph{separately}. 
The first cost function ${\rm Cost}_{\rm raw}$ is built out of \emph{raw data}, inserting exactly the voltages 
$V^{\rm experimental}(j,I_i)$ from the ``adjacent'' data set.  The second one, 
${\rm Cost}_{\rm corrected}$, is built 
out of \emph{explicitly correcting} the experimentally measured voltages $V^{\rm exprimental}$, by imposing an ansatz on how the errors 
across the different experiments are related. This pre-conditioning of the measured data results in 
the minimum of the cost function ${\rm Cost}_{\rm corrected}$ being much lower than that of the raw cost function ${\rm Cost}_{\rm raw}$. Moreover the visual quality of the reconstruction is most often similar, but at times much better than that obtained from the raw cost function. 
\medskip

\emph{Reconstruction of the contact impedances $z_k$:} It is known that errors in the contact impedances can lead to error in the reconstructions. This effect can be alleviated by performing the voltage measurements at  electrodes through which no current is injected~\cite{four-point-EIT}. In this work however we did not follow this approach, since we chose to use the entire ``adjacent'' data subset. 
We chose  to estimate the unknown contact impedances $z_k$, $k=1,\ldots,16$, \emph{first} 
by  using the measurements done in a tank filled only with saline. The $z_k$'s that were determined thus were then used in all subsequent experiments.  
\medskip

%to illustrate the effect of measurement errors on our reconstruction, we 
%perform a few of our experiments by making a ``by hand'' correction to the measured data for this data 
%subset: In particular we adjust the matrix corresponding to the adjacent data subset by subtracting from 
%each row its average. In a few instances, the visual quality of the reconstruction improves 
%substantially. In particular see experiments in Figure~\ref{fig:avg-error-removed} below.  For most of 
%the experiments, however, the value of the new cost function at the found minimum is much less than that 
%for the older cost function, but the visual quality of the reconstruction is not very different than the 
%old one. 
%\medskip

\emph{Dependence of ``found solution'' on the initial guess:}  Given the relatively small number of steps that we are allowing our optimizer to take, 
the final result depends on the initial guess; this is true for both the cost functions we minimize--the ``raw'' one and the ``corrected'' one. 
(This could of course also be explained by the cost function having many local minima--we are not sure if this is the case). 
For completeness, we include an example of the same problem solved twice: Once with a ``worse'' initial guess and 
once with a ``better'' one. (In both cases,   ``worse'' or ``better'' is in a visual sense). 
We find that when  we start with a better initial guess then (for the small number of steps we allow the optimizer to take), the reconstruction is  visually better. See Figure~\ref{fig:reco_perfect_guess}.
Crucially however, better visual fit at the end of the optimization also corresponds to better final ``scores'' for the cost functions. This gives substantial hope that more a tailored optimizer will be able to get even better results, without dramatically increasing the time needed (notwithstanding the issue of errors in the measurements). We also note that certain approaches to EIT that do not rely on optimization~\cite{BH-factorization-eit2,BH-factorization-eit3,Isaacson1,Isaacson2,d-bar1,d-bar2,BH-factorization-eit1,Kirsch-factorization1,Kirsch-factorization2,Nachman} could be used to provide a reasonable  first approximation for the solution, which could then be used as an initial guess for an inverse solver like the one we describe here. 
\medskip

\subsection{Figures of some reconstructions}\label{sec:real-data-reconstructions}

We present outcomes of our inverse solver applied to many of the ``phantoms'' considered in the FIPS
 experimental data.\footnote{ We remark that other inverse solvers have been tried on the same data set; the reader can see some reconstructions obtained by other inverse solvers in \cite{Experimental,OpenData}. }

The pictures we present below are the result of minimizing
\emph{two} cost functions: 
The first is the ``raw'' Cost function, which is defined in 
\eqref{cost final}, as well as a ``corrected'' Cost function, which we now describe.\footnote{The idea 
of the corrected cost function is entirely due to Joshua Fadelle, 
 as part of 
an undergraduate summer research experience at the University of Toronto; his work is 
co-funded by a Mathematics Department Research Award. } 

The ``corrected'' Cost function is given by the same formula \eqref{cost final}, by replacing the experimental 
voltages ${\rm V}^{\rm experimental}(j,I_i)$ by a corrected version 
${\rm V}^{\rm experimental}_{\rm corrected}(j,I_i)$, so 

\begin{equation}\label{cost corrected}
{\rm Cost}_{\rm corrected}[(\vec{B}_{\rm shapes},\vec{\sigma})]:= 10^6\cdot\sum_{i=1}^{16}\sum_{j=1}^{15}[
{\rm V}(j, I_i, \vec{B}_{\rm shapes}, \vec{\sigma})
-{\rm V}^{\rm experimental}_{\rm corrected}(j, I_i)]^2. 
\end{equation}
Here the ``corrected'' voltages are obtained by an ansatz on the errors in the measurements. 
We assume that for all the experiments, the current $I_i$ is measured correctly, yet at each pair of electrodes, given by $j\in \{1,\dots, 16\}$ the voltage $V^{\rm exprimental}(j,I_i)$
are measured with an error, which to high approximation is \emph{the same}, regardless of the 
 location of bodies. In other words, we assume the error is of the form $E(j,I_i)$. 
 We then identify the error $E(j,I_i)$ using the empty tank data. We identify the conductivity (and constant impedances) and find what the voltages \emph{should} have been for the found conductivity and contact impedances. We then define $E(j,I_i)$ to the difference between the experimentally measured 
 voltages and the values they should have had. 

 In particular we define
  \[
V^{\rm experimental}_{\rm corrected}(j,I_i):= V^{\rm experimental}(j,I_i)-E(j,I_i),
  \]
and define the function ${\rm Cost}_{\rm corrected}[(\vec{B}_{\rm shapes},\vec{\sigma})]$ 
using \eqref{cost corrected}.

We will present the results of minimizing \emph{both} cost functions, which we now denote by 
${\rm Cost}_{\rm raw}$ and 
${\rm Cost}_{\rm corrected}$ to make their distinction clear. 
Let us denote by 
\[
\vec{B}_{\rm found-shapes}^{\rm raw},~\vec{\sigma}_{\rm found-shapes}^{\rm raw}
\]
the results of the minimization for the ``raw'' cost function ${\rm Cost}_{\rm raw}$. We denote by 
\[
\vec{B}_{\rm found-shapes}^{\rm corrected},~\vec{\sigma}_{\rm found-shapes}^{\rm corrected}
\]
the results of the minimization for the ``corrected'' cost function ${\rm Cost}_{\rm corrected}$. 
The pictures of the found solutions for these two cost functions are displayed in figures 
\ref{fig:reco_1.1-1.2-1.3-1.4}--\ref{fig:reco_perfect_guess}. 
We next summarize what appears in those figures. 
\medskip

Recall that 
the minimization of either cost function  
that we perform requires an ``initial guess'' concerning the location of the bodies and 
their conductivities.

In the \emph{left} column 
we display  the \emph{initial guess} for $\vec{B}_{\rm shapes}$. (The blue curves superimposed on the pictures of the ``phantoms''). 
This initial guess is \emph{the same} for the minimization of both the ``raw'' and the ``corrected'' 
cost functions. 

In the \emph{middle} column 
we display the 
``found'' bodies  
$\vec{B}_{\rm found-shapes}^{\rm raw}$ (the red curves)
\emph{for the ``raw'' cost function ${\rm Cost}_{\rm raw}$, \eqref{cost final}}; i.e. the bodies corresponding to the obtained \emph{minimizer} 
\begin{equation}\label{found raw}
\vec{B}^{\rm raw}_{\rm found}= \big{(}\vec{B}^{\rm raw}_{\rm found-shapes},~\vec{\sigma}^{\rm raw}_{\rm found}\big{)}
\end{equation}
of the 
cost function \eqref{cost final}.

In the \emph{right} column we  display the 
``found'' bodies  
$\vec{B}_{\rm found-shapes}$ (the red curves)
\emph{for the ``corrected'' cost function ${\rm Cost}_{\rm corrected}$ \eqref{cost corrected}},
\begin{equation}\label{found corrected}
\vec{B}^{\rm corrected}_{\rm found}= \big{(}\vec{B}^{\rm corrected}_{\rm found-shapes},~\vec{\sigma}^{\rm corrected}_{\rm found}\big{)}.
\end{equation}

With labels we also note two parameters that relate to the optimization. 
 We list the value of the cost function \eqref{cost final}, \eqref{cost corrected} 
 at the \emph{found solutions} \eqref{found raw}, \eqref{found corrected} respectively; 
 we call this value the ``Score'' obtained by our optimizer. 
 For each experiment in the middle column  
 \[
 {\rm Score}= {\rm Cost}_{\rm raw}[(\vec{B}^{\rm raw}_{\rm found-shapes},~\vec{\sigma}^{\rm raw}_{\rm found})]
 \]
 and in the right-most column
 \[
 {\rm Score}= {\rm Cost}_{\rm corrected}[(\vec{B}^{\rm corrected}_{\rm found-shapes},~\vec{\sigma}^{\rm corrected}_{\rm found})]. 
 \]
 We also list the number iterations $i$ that  the optimizer took.  
 \medskip

  %In the first list of reconstructions, Figures~\ref{fig:reco_1.1-1.2-1.3}-\ref{fig:reco_5.2} the cost function is constructed 
 % exclusively from the \emph{raw} ``adjacent data'' in the FIPS data set. 

The initial guesses in the left columns of Figures~\ref{fig:reco_1.1-1.2-1.3-1.4}-\ref{fig:reco_5.1-5.2}  are chosen essentially at random, trying to make sure they are reasonably far from the (visually) ``correct'' locations. The reader will note that \emph{a few} of the phantoms are not ``found'' very well, notably in the last two examples in Figure~\ref{fig:reco_2.5-2.6-3.1-3.2} or the first two in Figure~\ref{fig:reco_4.1-4.2-4.3-4.4} the reconstruction from ``raw'' data is not as successful.

%The first two such  are phantoms 4.1, 5.2; for those we also present the outcome of minimizing a 
%\emph{different} cost function: For 4.1, 5.2 
%we define a new cost function, by  using instead of the 
%``raw'' data from the ``adjacent'' data subset, a \emph{modification} of this data, where we  modify
%each row in the data matrix by subtracting its average. 
%The (visually better) results for these two experiments are in Figure~\ref{fig:avg-error-removed}. 
%\medskip

We next present some evidence that in the examples where the ``raw data'' reconstruction was not very good, it is likely that the optimizer found only a local minimum. Our evidence for this is based on performing a new optimization for the ``poor'' reconstruction of the last setting in Figure~\ref{fig:reco_2.5-2.6-3.1-3.2} (for the same cost functions ${\rm Cost}_{\rm raw}$ and ${\rm Cost}_{\rm corrected}$) but changing the initial guess. We note that the found solutions are better, both visually, and (crucially)  in terms of the lower value of the Cost functions at the found solution. This is performed in Figure~\ref{fig:reco_perfect_guess}.

\subsection{Conclusion} 
We introduced a novel boundary integral equation based solver for the electrostatic potential equation \eqref{thePDE}, 
which we showed to be both fast and accurate. We used this forward solver to construct  a bare-bones inverse solver for the EIT inverse problem, which we then tested on the FIPS open data set. Applying the EIT inverse solver on the raw data yields very satisfactory reconstructions and using it on suitably pre-processed data can either maintain the quality of the reconstruction or even improve it markedly. 

The EIT method we introduced lends itself to further generalizations. The inverse solver could be improved by a more tailored choice of cost function. One could utilize faster optimization algorithms or algorithms that seek to also account for the error in the measurements as part for the inverse solver. Of course, one could also seek to extend the methods to three dimensions to capture more potential applications. These challenges will be undertaken in future research.
  
%We also present the output of minimizing the same ``average-adjusted'' cost function, but in addition improving the initial guess to be (visually) closer to the correct solution. 
%This is done for phanto 5.2. The results of this are in Figure~\ref{fig:reco_perfect_guess}. We note that the ``scores'' obtained are substantially better than the scores with the ``visually worse'' initial guess. We believe this should be interpreted as  the solver previously having gotten stuck in a local minimum, or simply terminated due to the maximum number of allowed steps being reached. The found images Figure~\ref{fig:reco_perfect_guess} are thus perhaps closer to the global minimum. 
\medskip

\newcommand{\figwidth}{0.85}
\begin{figure}[p]
    \centering
    \includegraphics[width=\figwidth\textwidth]{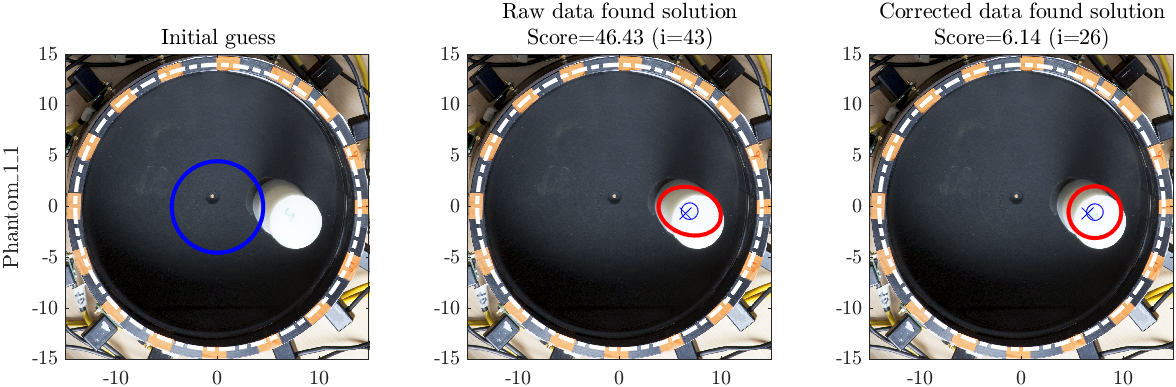}
\smallskip
    \includegraphics[width=\figwidth\textwidth]{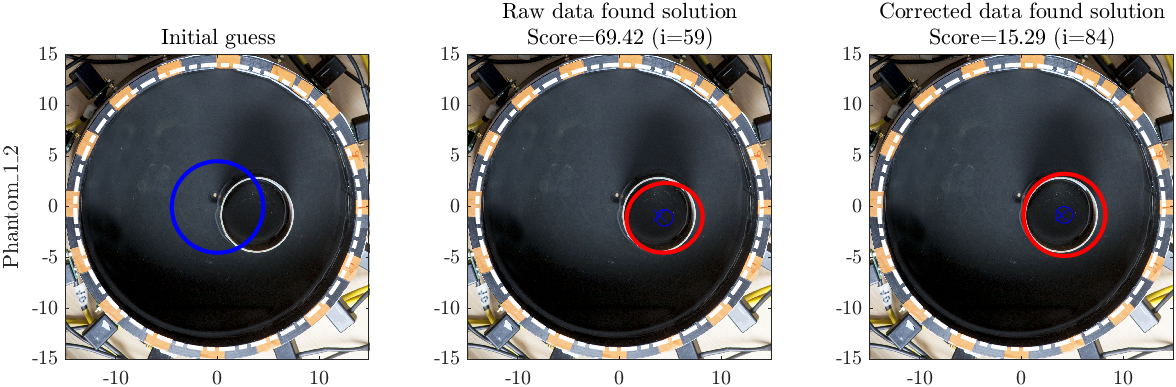}
\smallskip
    \includegraphics[width=\figwidth\textwidth]{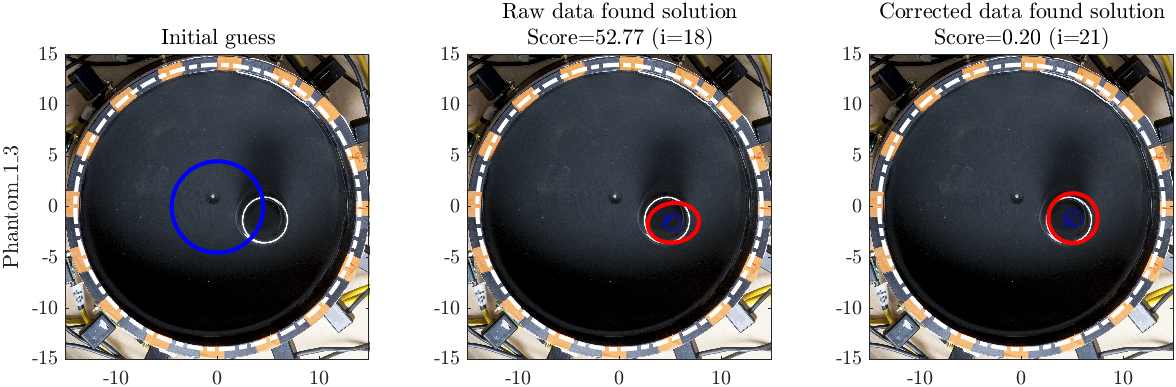}
    \smallskip
    \includegraphics[width=\figwidth\textwidth]{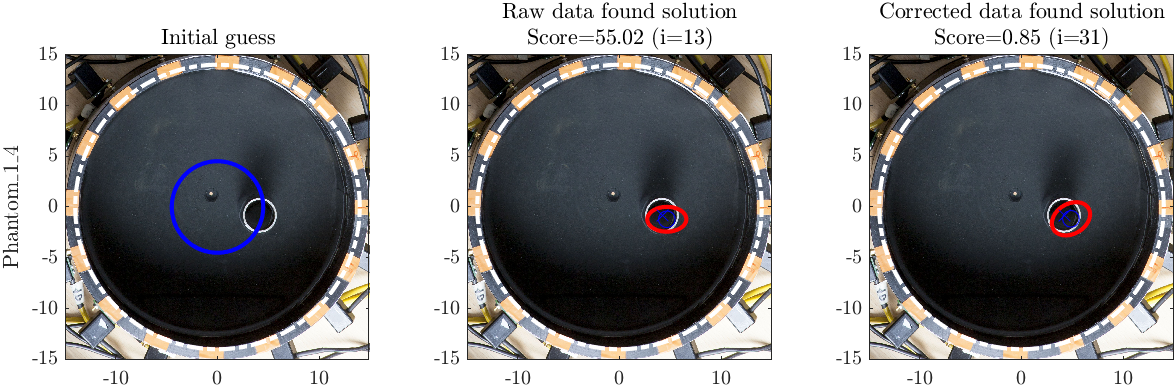}
    \caption{The initial guess on the left; the found solution from raw data in the middle; the found solution from corrected data on the right. The score values of the cost functions ${\rm Cost}_{\rm raw}$ and ${\rm Cost}_{\rm corrected}$ at the found solutions are displayed above each photo. The visually estimated and the numerically found centers of masses for the bodies are marked by a blue cross and a circle, respectively.\label{fig:reco_1.1-1.2-1.3-1.4}}
\end{figure}

\begin{figure}[p]
    \centering
    \includegraphics[width=\figwidth\textwidth]{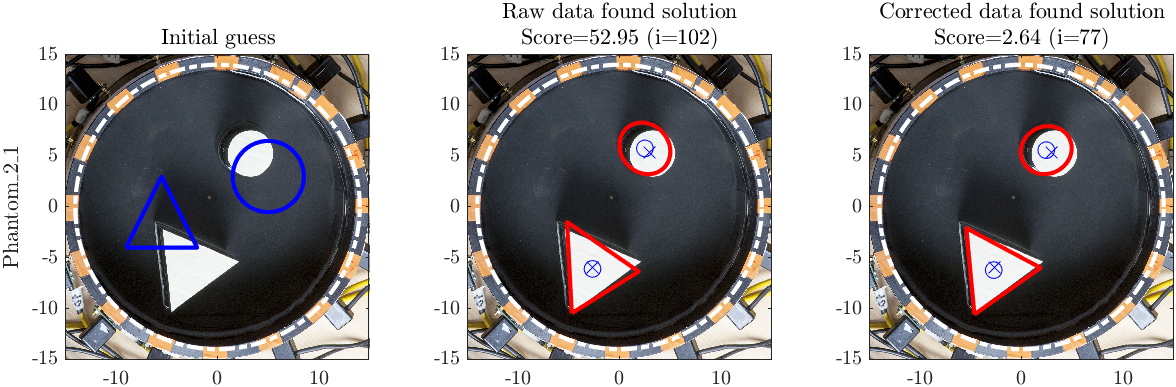}
\smallskip
    \includegraphics[width=\figwidth\textwidth]{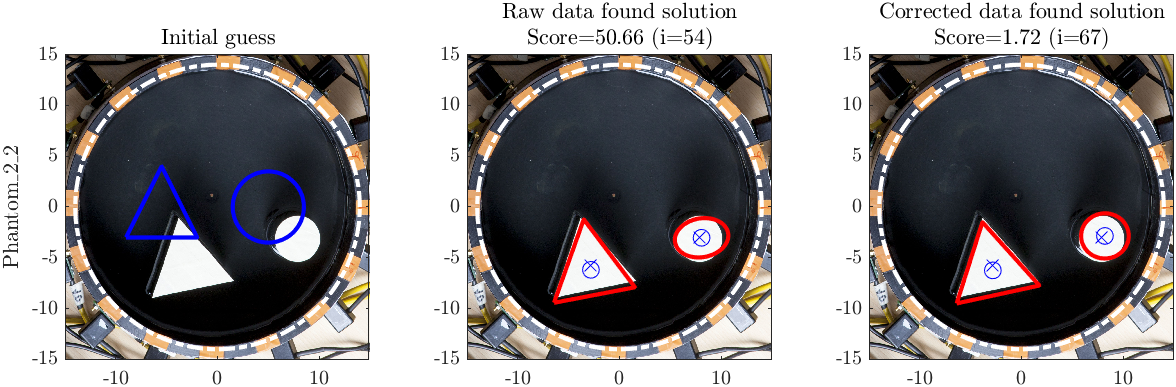}
\smallskip
    \includegraphics[width=\figwidth\textwidth]{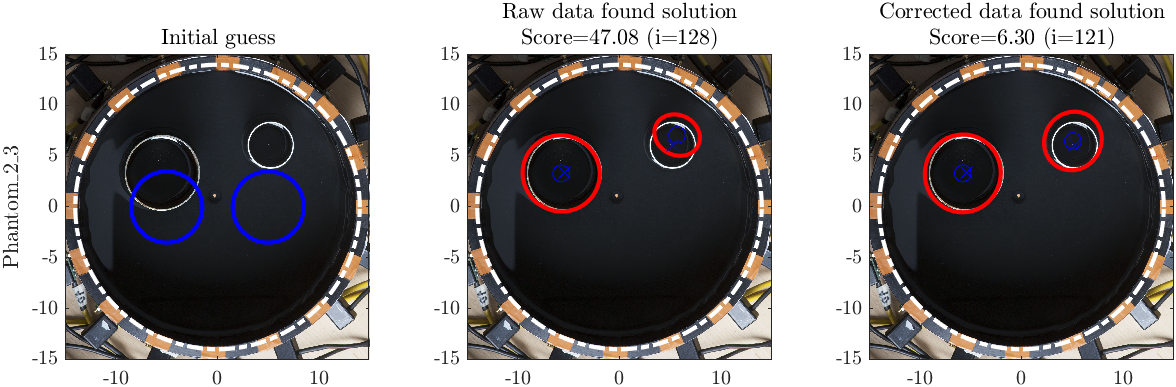}
    \smallskip
    \includegraphics[width=\figwidth\textwidth]{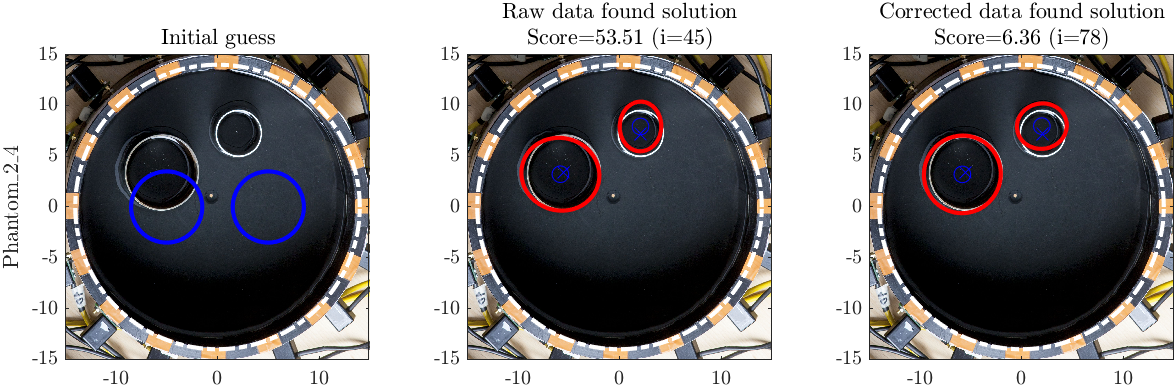}
    \caption{The initial guess on the left, the found solution from raw data in the middle, the found solution from corrected data on the right. The score values of the cost functions ${\rm Cost}_{\rm raw}$ and ${\rm Cost}_{\rm corrected}$ at the found solutions displayed above each figure. The visually estimated and the numerically found centers of masses for the bodies are marked by a blue cross and a circle, respectively.\label{fig:reco_2.1-2.2-2.3-2.4}}
\end{figure}

\begin{figure}[p]
    \centering
    \includegraphics[width=\figwidth\textwidth]{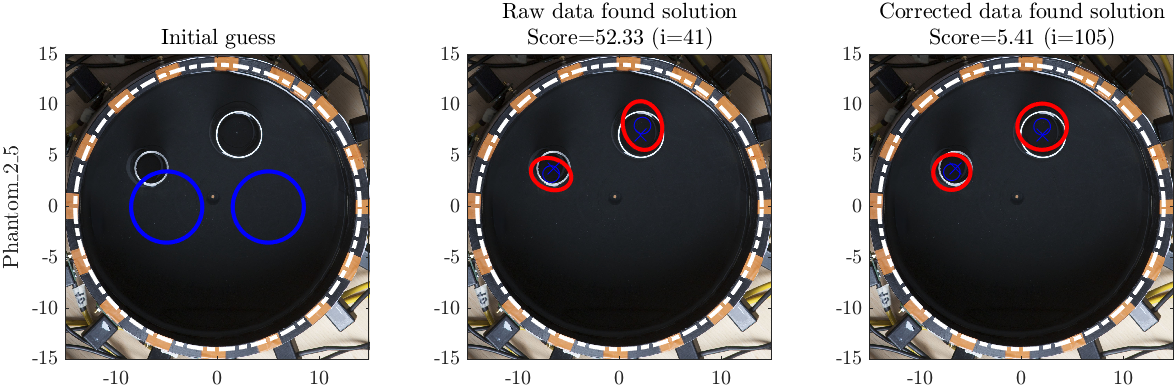}
\smallskip
    \includegraphics[width=\figwidth\textwidth]{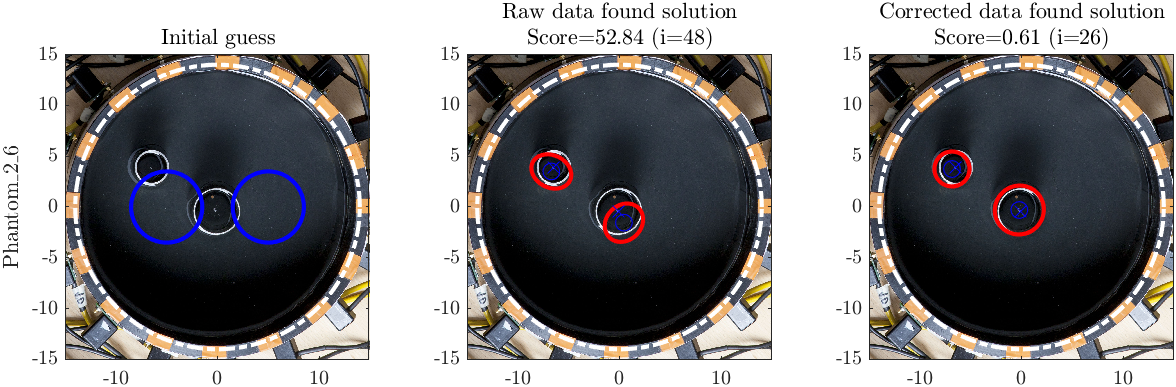}
\smallskip
    \includegraphics[width=\figwidth\textwidth]{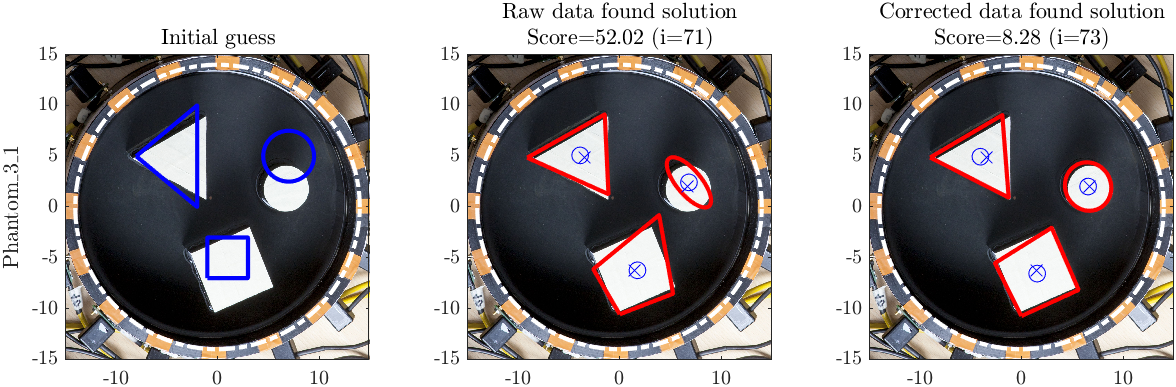}
    \smallskip
    \includegraphics[width=\figwidth\textwidth]{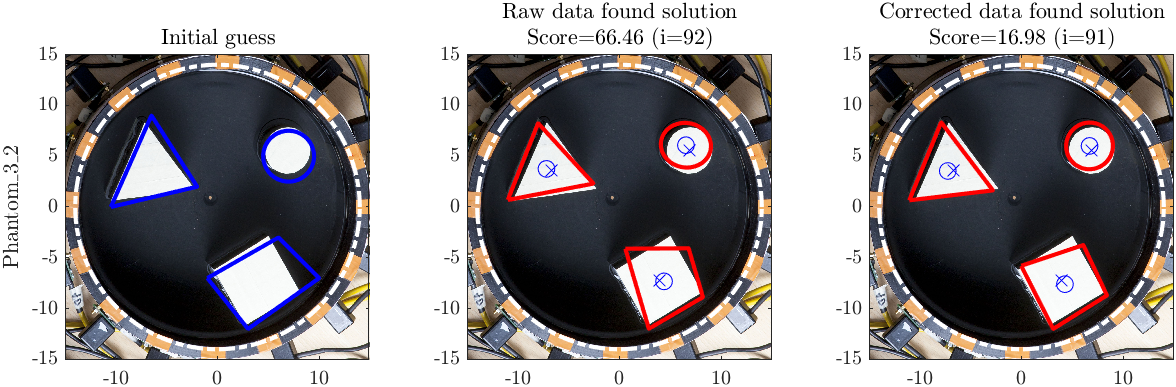}
    \caption{The initial guess on the left; the found solution from raw data in the middle; the found solution from corrected data on the right. The score values of the cost functions ${\rm Cost}_{\rm raw}$ and ${\rm Cost}_{\rm corrected}$ at the found solutions are displayed above each photo. The visually estimated and the numerically found centers of masses for the bodies are marked by a blue cross and a circle, respectively.\label{fig:reco_2.5-2.6-3.1-3.2}}
\end{figure}

\begin{figure}[p]
    \centering
    \includegraphics[width=\figwidth\textwidth]{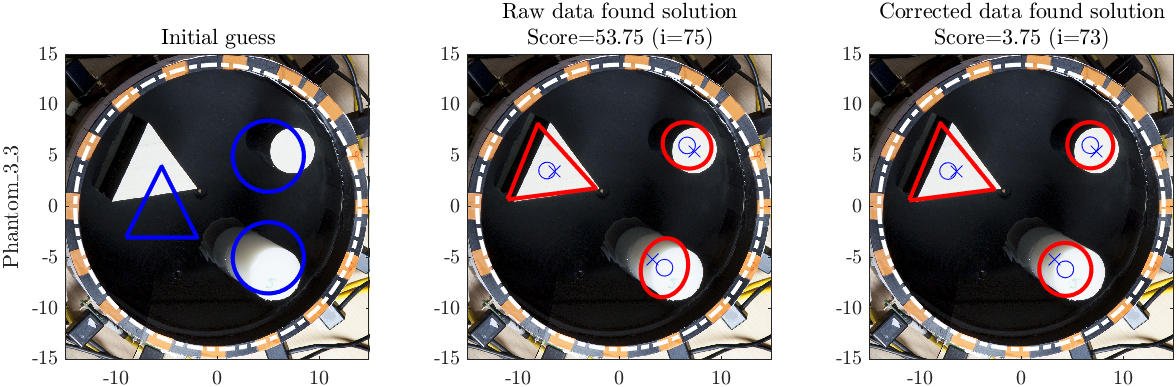}
\smallskip
    \includegraphics[width=\figwidth\textwidth]{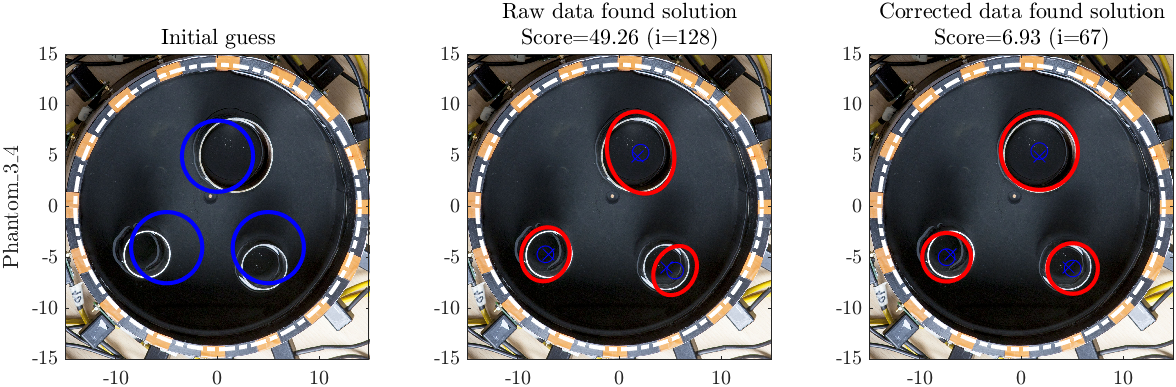}
\smallskip
    \includegraphics[width=\figwidth\textwidth]{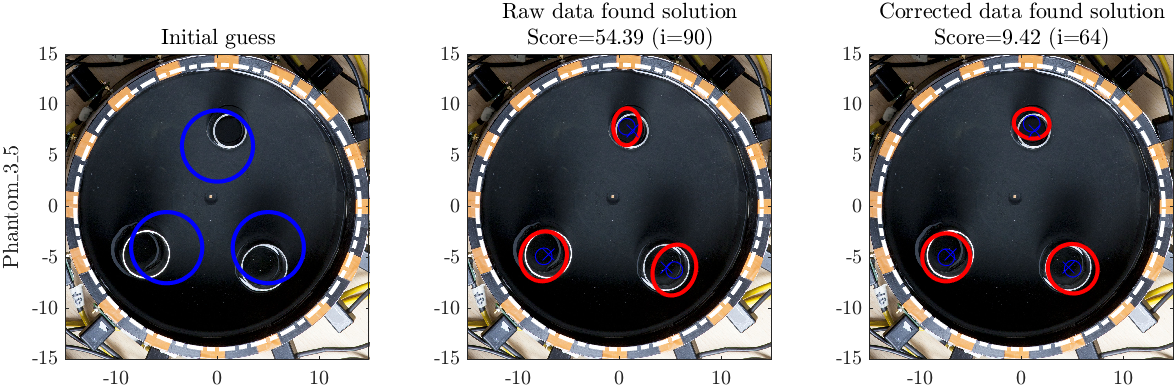}
    \smallskip
    \includegraphics[width=\figwidth\textwidth]{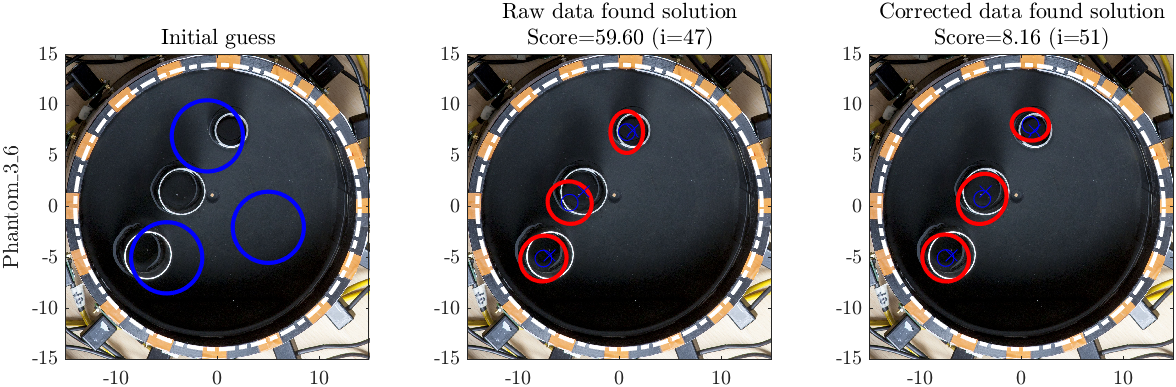}
    \caption{The initial guess on the left; the found solution from raw data in the middle; the found solution from corrected data on the right. The score values of the cost functions ${\rm Cost}_{\rm raw}$ and ${\rm Cost}_{\rm corrected}$ at the found solutions are displayed above each photo. The visually estimated and the numerically found centers of masses for the bodies are marked by a blue cross and a circle, respectively.\label{fig:reco_3.3-3.4-3.5-3.6}}
\end{figure}

\begin{figure}[p]
    \centering
    \includegraphics[width=\figwidth\textwidth]{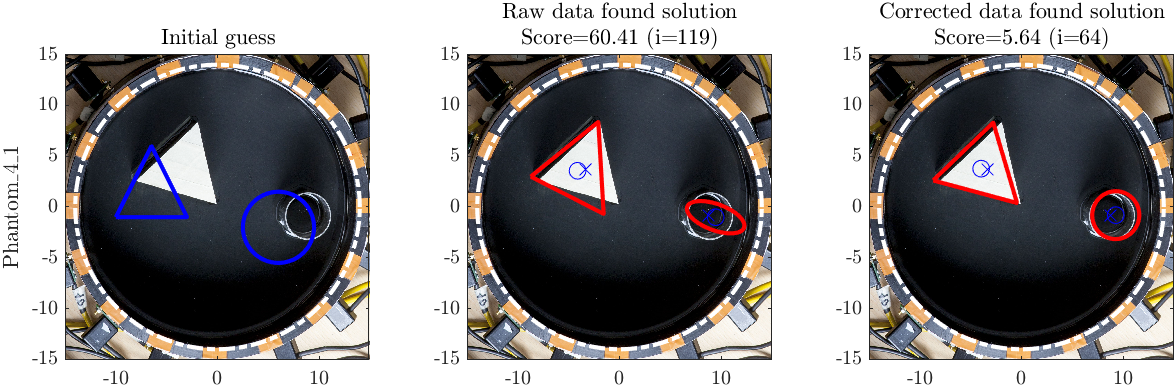}
\smallskip
    \includegraphics[width=\figwidth\textwidth]{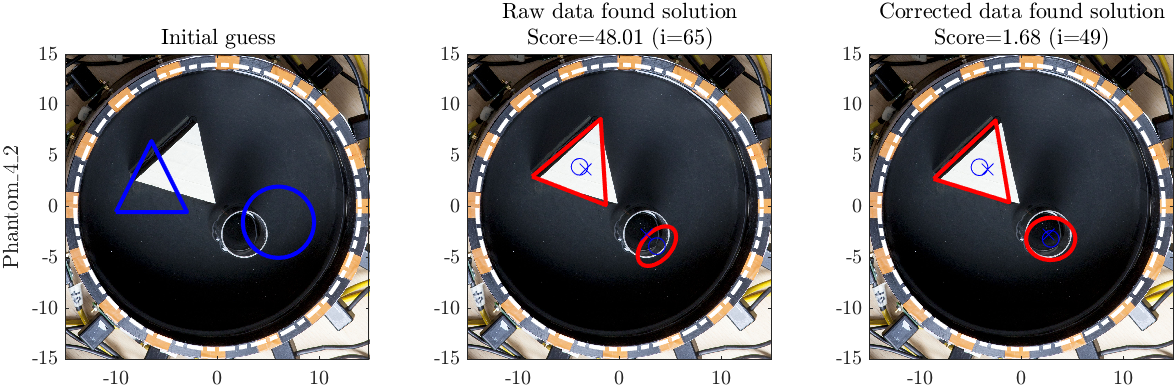}
\smallskip
    \includegraphics[width=\figwidth\textwidth]{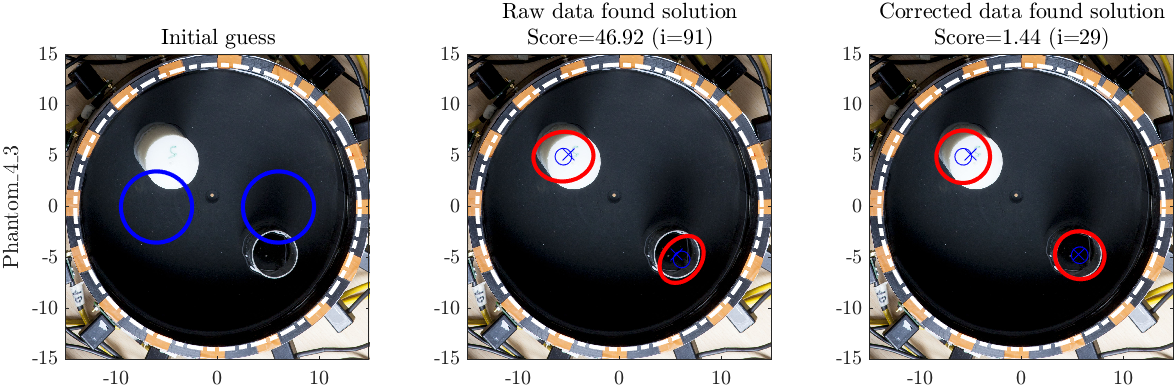}
    \smallskip
    \includegraphics[width=\figwidth\textwidth]{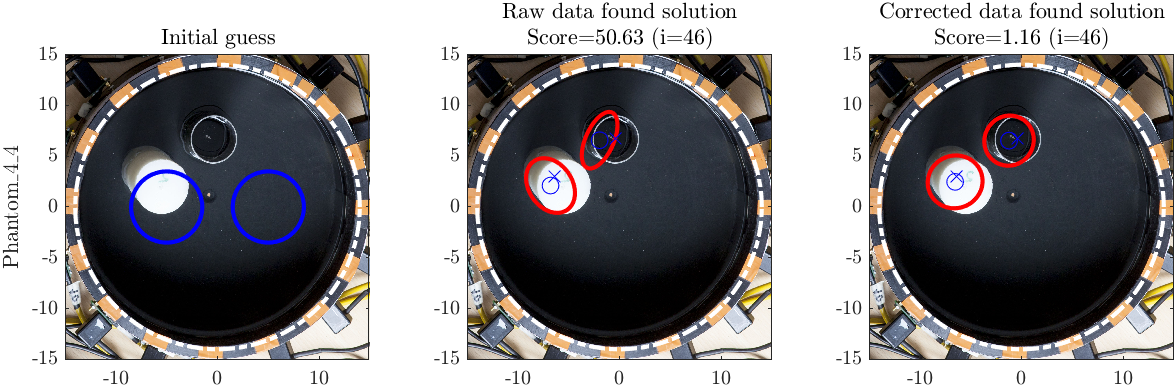}
    \caption{The initial guess on the left; the found solution from raw data in the middle; the found solution from corrected data on the right. The score values of the cost functions ${\rm Cost}_{\rm raw}$ and ${\rm Cost}_{\rm corrected}$ at the found solutions are displayed above each photo. The visually estimated and the numerically found centers of masses for the bodies are marked by a blue cross and a circle, respectively.\label{fig:reco_4.1-4.2-4.3-4.4}}
\end{figure}

\begin{figure}[p]
    \centering
    \includegraphics[width=\figwidth\textwidth]{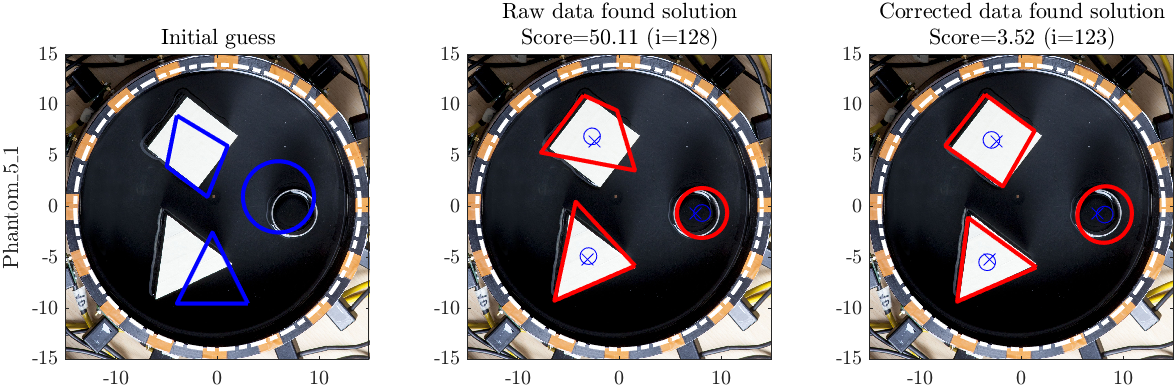}
\smallskip
    \includegraphics[width=\figwidth\textwidth]{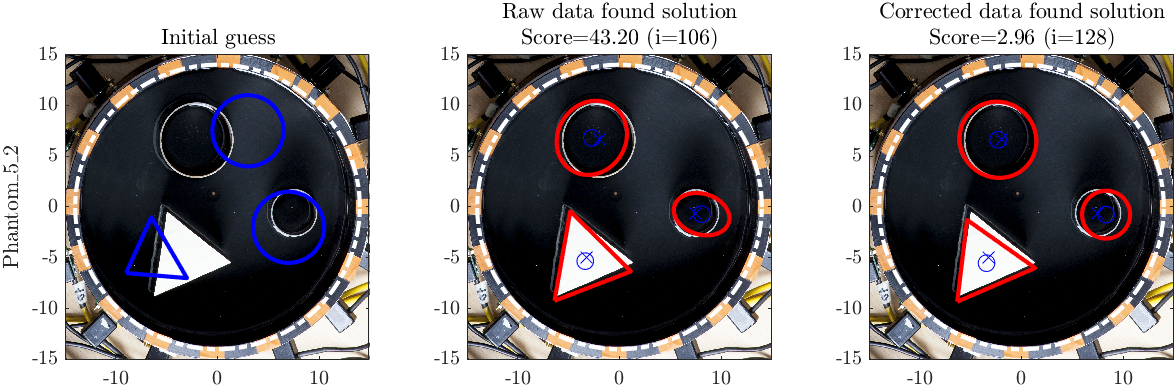}
    \caption{The initial guess on the left; the found solution from raw data in the middle; the found solution from corrected data on the right. The score values of the cost functions ${\rm Cost}_{\rm raw}$ and ${\rm Cost}_{\rm corrected}$ at the found solutions are displayed above each photo. The visually estimated and the numerically found centers of masses for the bodies are marked by a blue cross and a circle, respectively.\label{fig:reco_5.1-5.2}}
\end{figure}

\begin{figure}[p]
    \centering
    \includegraphics[width=\figwidth\textwidth]{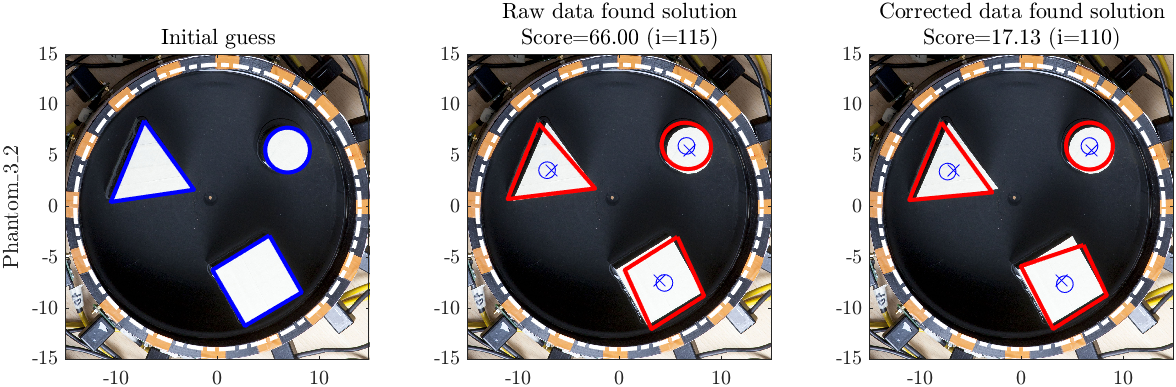}
    \caption{The initial guess on the left; the found solution from raw data in the middle; the found solution from corrected data on the right. Visually perfect initial guess can yield better reconstructions and scores, compared to the last example of Figure~\ref{fig:reco_2.5-2.6-3.1-3.2}. \label{fig:reco_perfect_guess}}
\end{figure}

\subsection{Acknowledgements}
The authors would like to thank Kirill Serkh for many helpful discussions. We also thank Kyle Bower for his help on implementing the overlap conditions of interfaces in the BIEM, done for his PhD thesis~\cite{Bower-thesis}, see also \cite{Bower}. We are also grateful to Joshua Fadelle for his insight for the corrected cost function. This work was supported by the Natural Sciences and Engineering Research Council of Canada (NSERC): RGPIN-2020-01113 and RGPIN-2019-06946.

\end{document}